# Mordukhovich derivatives of the set-valued metric projection operator in general Banach spaces

Jinlu Li

Department of Mathematics
Shawnee State University
Portsmouth, Ohio 45662
USA

**Abstract**

In this paper, we investigate the properties and the precise solutions of the Mordukhovich derivatives of the set-valued metric projection operator onto some closed balls in some general Banach spaces. In the Banach space $c$, we find the properties of Mordukhovich derivatives of the set-valued metric projection operator onto the closed subspace $c_0$. We show that the metric projection from $C[0, 1]$ to polynormal with degree less than or equal to $n$ is a single-valued mapping. We investigate its Mordukhovich derivatives and Gâteaux directional derivatives.

1. ## Introduction

Let $(X, \|\cdot\|)$ be a Banach space with topological dual space $(X^*, \|\cdot\|_*)$. Let $C$ be a nonempty closed and convex subset of $X$ and let $\theta$ be the origin of $X$. Let $P_C$ denote the (standard) metric projection operator from $X$ to $C$. It is well-known that, in the operator theory, the metric projection operator is one of the most important operators, which has been studied by many authors with a long history (see [1, 2, 3, 6, 7, 10, 13, 18, 23, 27, 28]). The metric projection operator has been extensively applied to many fields in nonlinear analysis such as, optimization theory, approximation theory, fixed point theory, variational analysis, equilibrium theory, and so forth (see [1, 2, 6, 26, 27, 28]).

If the considered Banach space is uniformly convex and uniformly smooth, then the metric projection operator $P_C: X \to C$ is a single-valued mapping. In this case, we can consider the continuity of $P_C$, which is one of the most important research topics in operator theory. One more step further, in addition to the continuity of $P_C$, the smoothness of $P_C$ has attracted a lot of attention in the fields of nonlinear analysis. Traditionally, the smoothness of $P_C$ is described by some differentiability. For example, if the underlying spaces are Hilbert spaces, some differentiability of $P_C$ has been introduced and studied in [3, 6, 7, 10, 11, 14, 15, 19], and if the underlying spaces are Banach spaces and normed linear spaces in [4, 9, 14, 17, 26, 27].

We believe that the most popular and useful concept of differentiability of $P_C$ is the Fréchet differentiability. We recall the concepts of Gâteaux directional differentiability and the Fréchet differentiability of the metric projection operator $P_C: X \to C$, in which the underlying space is a uniformly convex and uniformly smooth Banach space (see Definition 4.1 in [12] and Definition 1.13 in [19]). For $\bar{x} \in X$ and $v \in X$ with $v \neq \theta$, if the following limit exists,

$$P'_C(\bar{x})(v) := \lim_{t \downarrow 0} \frac{P_C(\bar{x}+tv) - P_C(\bar{x})}{t},$$

then, $P_C$ is said to be Gâteaux directionally differentiable at point $\bar{x}$ along direction $v$ and $P'_C(\bar{x})(v)$ is called the Gâteaux directional derivative of $P_C$ at point $\bar{x}$ along direction $v$.

If there is a linear and continuous mapping $\nabla P_C(\bar{x}): X \to X$ such that

$$\lim_{x \to \bar{x}} \frac{P_C(x) - P_C(\bar{x}) - \nabla P_C(\bar{x})(x - \bar{x})}{\|x - \bar{x}\|} = \theta,$$

then $P_C$ is said to be Fréchet differentiable at $\bar{x}$ and $\nabla P_C(\bar{x})$ is called the Fréchet derivative of $P_C$ at $\bar{x}$. The above two types of differentiability have the following inclusion properties (see [14] for mor details): for any given point $\bar{x} \in X$,

$P_C$ is Fréchet differentiable at $\bar{x}$ $\implies$ $P_C$ is Gâteaux directionally differentiable at $\bar{x}$.

However, in general, we have

$P_C$ is Gâteaux directionally differentiable at $\bar{x}$ $\not\Rightarrow$ $P_C$ is Fréchet differentiable at this point $\bar{x}$.

We notice that the concepts of Gâteaux directional differentiability and Fréchet differentiability of the metric projection operator inherited the ordinary definitions of differentiation in functional analysis, which is considered as the traditional differentiability of the metric projection $P_C$.

However, if the underlying space is a general Banach space, then, the metric projection operator may be a set-valued mapping. Consequently, the continuity and the ordinary differentiability of the metric projection $P_C$ cannot be traditionally defined as the above definition of Fréchet differentiability.

In nonlinear analysis, the set-valued mappings show increasingly important roles, such as, in variational analysis, game theory, equilibrium theory, economic theory, and so forth. In [19, 20, 21], Mordukhovich successfully developed the generalized and fresh theory of differentiability of set-valued mappings (Which include single-valued mappings as special cases) in Banach spaces, which is called generalized differentiation in [19] (see Definitions 1.13 and 1.32 in Chapter 1 in [19]) and we call it Mordukhovich differentiation. The Mordukhovich derivative construction introduced in [19] has become a fundamental method of modern analysis theory. It has widely influenced several branches of mathematics such as operator theory, optimization theory, approximation theory, control theory, equilibrium theory, and so forth.

In particular, if the underlying Banach space is uniformly convex and uniformly smooth, then, the metric projection operator $P_C$ becomes a single-valued mapping. The Mordukhovich coderivative of $P_C$ at a point $(\bar{x}, P_C(\bar{x}))$ is denoted by $\widehat{D}^* P_C(\bar{x})$ that is defined to be a set-valued mapping from $X^*$ to $X^*$ such that, for any $y^* \in X^*$,

$$\widehat{D}^* P_C(\bar{x})(y^*) = \left\{ x^* \in X^* : \limsup_{u \to \bar{x}} \frac{\langle x^*, u - \bar{x} \rangle - \langle y^*, P_C(u) - P_C(\bar{x}) \rangle}{\|u - \bar{x}\| + \|P_C(u) - P_C(\bar{x})\|} \leq 0 \right\}.$$

Furthermore, if the considered mapping is a single-valued mapping (which is considered as a special case of set-valued mapping), then, Theorem 1.38 in [19] shows that there are very close connections between the Fréchet derivative $\nabla P_C(\bar{x})$ and the Mordukhovich derivative at this given point $\bar{x}$. By the results of Theorem 1.38 in [19], we can consider that the Mordukhovich differentiation generalized the Fréchet differentiation from single-valued mappings to set-valued mappings in Banach spaces.

Very recently, the present author investigated the properties of Mordukhovich derivatives of the metric projection operator in Hilbert spaces in [16], and in uniformly convex and uniformly smooth Banach spaces in [18], respectively. In section 3 of this paper, we investigate the properties and the precise solutions of the Mordukhovich derivatives of the set-valued metric projection operator onto some closed balls in some general Banach spaces. In section 4, we consider the Banach space $c$ and we find the properties of Mordukhovich derivatives of the set-valued metric projection operator onto the closed subspace $c_0$. In section 5, we study the metric projection operator from real Banach space $C[0, 1]$ to polynormal with degree less than or equal to $n$. We prove that it is a single-valued mapping and investigate its Mordukhovich derivatives and Gâteaux directional derivatives.

## 2. Preliminaries

### 2.1. The normalized duality mapping and its inverse in general Banach spaces

Let $(X, \|\cdot\|)$ be a real Banach space with topological dual space $(X^*, \|\cdot\|_*)$. Let $\langle \cdot, \cdot \rangle$ denote the real canonical pairing between $X^*$ and $X$, and $\theta, \theta^*$ denote the origins in $X$ and $X^*$, respectively. Let $\mathbb{B}$ and $\mathbb{B}^*$ denote the unit closed balls in $X$ and $X^*$, respectively. It follows that, for any $r > 0$, $r\mathbb{B}$ and $r\mathbb{B}^*$ are closed balls with radius $r$ and center origins in $X$ and $X^*$, respectively. Let $\mathbb{S}$ be the unit sphere in $X$. Then, $r\mathbb{S}$ is the sphere in $X$ with radius $r$ and center $\theta$. For any $c \in X$ and $r > 0$, let $\mathbb{B}(c, r)$ denote the closed ball in $X$ with radius $r$ and center $c$. It follows that $\mathbb{B}(\theta, 1) = \mathbb{B}$ and $\mathbb{B}(\theta, r) = r\mathbb{B}$. The identity mappings on $X$ and $X^*$ are respectively denoted by $I_X$ and $I_{X^*}$.

The normalized duality mapping $J: X \to 2^{X^*} \setminus \{\emptyset\}$ is defined by

$$J(x) = \{jx \in X^*: \langle jx, x \rangle = \|jx\|_* \|x\| = \|x\|^2 = \|jx\|_*^2\}, \text{ for any } x \in X.$$

The normalized duality mapping on $X^*$ is denoted by $J^*$. In general, when the considered Banach space is not reflexive, for any $x^* \in X^*$, the value $J^*(x^*)$ contains a subset in $X$. For $x, y \in X$, if $J(x) \cap J(x) \neq \emptyset$, then $x$ and $y$ are said to be generalized identical. For $x \in X$, the set of all its generalized identical points is denoted by $\mathfrak{I}(x)$. That is,

$$\mathfrak{I}(x) = \{y \in X: J(x) \cap J(x) \neq \emptyset\}.$$

Some properties of the normalized duality mapping and generalized identity in general Banach spaces are listed in the appendix of this paper.

### 2.2. The variational properties of the metric projection on general Banach spaces

Let $C$ be a nonempty closed and convex subset of Banach space $X$. Let $P_C: X \to 2^C$ denote the (standard) metric projection operator, which is a set-valued mapping, in general. For any $x \in X$,

we have $P_C(x) \subseteq C$ such that

$$P_C(x) = \{u \in C: \|x - u\| \leq \|x - z\|, \text{ for all } z \in C\}.$$

However, when $P_C$ is restricted on $C$, then, it is single-valued such that

$$P_C(x) = x, \text{ for every } x \in C.$$

The metric projection operator $P_C$ has the following basic variational principle (see [11]). A proof of the following proposition is provided in the appendix of this paper.

**Proposition 2.1.** *Let X be a Banach space and C a nonempty closed and convex subset of X. For any $x \in X$ and $u \in C$, then, $u \in P_C(x)$, if and only if, there is $j(x-u) \in J(x-u)$ such that*

$$\langle j(x-u), u-y \rangle \geq 0, \quad \text{for all } y \in C.$$

### 2.3. The product spaces of Banach spaces

Let $X \times X$ and $X^* \times X^*$ be the Cartesian product spaces of $X$ and $X^*$, respectively. In this paper, the norms on $X \times X$ and $X^* \times X^*$ are respectively denoted by $\|\cdot\|_A$ and $\|\cdot\|_{A^*}$. We will define the real canonical pairing between $X^* \times X^*$ and $X \times X$ such that the normalized duality mapping on $X \times X$ will be defined with respect to the norms $\|\cdot\|_A$ and $\|\cdot\|_{A^*}$ and the canonical pairing. See [24] for more details about product spaces.

In this paper, the norm $\|\cdot\|_A$ on $X \times X$ is defined by

$$\|(x,y)\|_A = \sqrt{\|x\|^2 + \|y\|^2}, \text{ for all } (x,y) \in X \times X. \tag{2.1}$$

Similarly, the norm $\|\cdot\|_{A^*}$ on $X^* \times X^*$ is defined by

$$\|(x^*,y^*)\|_{A^*} = \sqrt{\|x^*\|_*^2 + \|x^*\|_*^2}, \text{ for all } (x^*,y^*) \in X^* \times X^*. \tag{2.2}$$

Let $\langle \cdot, \cdot \rangle_A$ denote the real canonical pairing between $X^* \times X^*$ and $X \times X$, which is defined, for any $(x,y) \in X \times X$ and $(x^*,y^*) \in X^* \times X^*$, by

$$\langle (x^*,y^*), (x,y) \rangle_A = \langle x^*, x \rangle + \langle y^*, y \rangle. \tag{2.3}$$

The Cartesian product space $X \times X$ with the norm $\|\cdot\|_A$ is also a Banach space with the dual space $X^* \times X^*$ and with the real canonical pairing $\langle \cdot, \cdot \rangle_A$ between $X^* \times X^*$ and $X \times X$.

Based on $J: X \to 2^{X^*}$, we define a mapping $\mathbb{J}: X \times X \to 2^{X^* \times X^*}$, for any $(x,y) \in X \times X$, by

$$\mathbb{J}(x,y) := \{(j(x), j(y)) \in X^* \times X^*: j(x) \in J(x) \text{ and } j(y) \in J(y)\}.$$

By (2.1) and (2.3), for any given $(x,y) \in X \times X$, for any $(j(x), j(y)) \in \mathbb{J}(x,y)$, we have

$$\langle (j(x), j(y)), (x,y) \rangle_A$$
$$= \langle j(x), x \rangle + \langle j(y), y \rangle$$

$$= \|x\|^2 + \|y\|^2$$
$$= \|(x,y)\|_A^2. \tag{2.4}$$

On the other hand, by (2.2) and (2.3), we have

$$\langle (j(x), j(y)), (x, y)\rangle_A$$
$$= \|x\|^2 + \|y\|^2$$
$$= \|jx\|_*^2 + \|jy\|_*^2$$
$$= \|(j(x), j(y))\|_{A^*}^2. \tag{2.5}$$

By (2.4) and (2.5), we obtain that, for any $(j(x), j(y)) \in \mathbb{J}(x, y)$, $(j(x), j(y))$ satisfies the conditions to be a normalized duality mapping from $X \times X$ to $X^* \times X^*$. This implies that the mapping $\mathbb{J}: X \times X \to 2^{X^* \times X^*}$ is the normalized duality mapping from $X \times X$ to $X^* \times X^*$.

### 2.4. Mordukhovich derivatives of the metric projection operator in Banach spaces

In [19], the concepts of Mordukhovich derivatives of set-valued mappings in Banach spaces are introduced. In this paper, we focus on investigating Mordukhovich derivatives of the metric projection operator in Banach spaces. Here, we rewrite the definitions of Mordukhovich derivatives given in [19] with respect to the set-valued mapping $P_C: X \rightrightarrows C$.

**Definition 1.1 in [19]** (generalized normals in Banach spaces). Let $\Omega$ be a nonempty subset of a Banach space $X$. For any $x \in X$, we define

$$\widehat{N}(x; \Omega) = \left\{ x^* \in X^* : \limsup_{\substack{u \xrightarrow{\Omega} x}} \frac{\langle x^*, u-x \rangle}{\|u-x\|} \leq 0 \right\}, \text{ for any } x \in \Omega. \tag{2.6}$$

The elements of (2.6) are called Fréchet normals and $\widehat{N}(x; \Omega)$ is the prenormal cone to $\Omega$ at $x$. We put $\widehat{N}(x; \Omega) = \emptyset$, for any $x \in X \setminus \Omega$. Then, based on the above definition, we recall the Mordukhovich derivatives of the set-valued mapping $P_C$ in Banach spaces. Note that, in [19], Mordukhovich derivatives are named by Fréchet coderivatives or precoderivatives.

**Definition 1.32 in [19].** Let $C$ be a nonempty closed and convex subset of Banach space $X$. For any $x \in X$ and $y \in P_C(x)$, the Mordukhovich derivatives of $P_C$ at $(x, y)$ is defined to be a set-valued mapping $\widehat{D}^* P_C(x, y): X^* \rightrightarrows X^*$ such that, for any given $y^* \in X^*$,

$$\widehat{D}^* P_C(x, y)(y^*) = \{x^* \in X^* : (x^*, -y^*) \in \widehat{N}((x,y); \text{gph} P_C)\}. \tag{2.7}$$

If $u \notin P_C(x)$, then, we define $\widehat{D}^* P_C(x, u)(y^*) = \emptyset$, for all $y^* \in X^*$. By Definition 1.1 in [19], we convert the Definition (2.7) to a limit for calculating the Mordukhovich derivative of $P_C$ at $(x, y)$ as follows. For any $y^* \in X^*$, we have

$$\widehat{D}^*P_C(x,y)(y^*) = \left\{ x^* \in X^*: \limsup_{\substack{(u,v) \to (x,y) \\ v \in P_C(u)}} \frac{\langle (x^*,-y^*),(u,v)-(x,y) \rangle_A}{\|(u,v)-(x,y)\|_A} \leq 0 \right\}$$

$$= \left\{ x^* \in X^*: \limsup_{\substack{(u,v) \to (x,y) \\ v \in P_C(u)}} \frac{\langle x^*, u-x \rangle - \langle y^*, v-y \rangle}{\sqrt{\|u-x\|^2 + \|v-y\|^2}} \leq 0 \right\}. \tag{2.8}$$

As mentioned in the proof of Proposition 12 in Chapter 1 of [19], the limit in (2.8) does not depend on equivalent norms on $X \times X$. Notice that

$$\frac{\sqrt{2}}{2}(\|x\| + \|y\|) \leq \sqrt{\|x\|^2 + \|y\|^2} \leq \|x\| + \|y\|, \text{ for all } (x,y) \in X \times X. \tag{2.9}$$

Hence, in the limit (2.8), we replace the norm on $X \times X$ by the equivalent norm in (2.9)

$$\|(x,y)\| = \|x\| + \|y\|, \text{ for all } (x,y) \in X \times X.$$

For any $x \in X$ and $y \in P_C(x)$, (2.8) is rewritten as

$$\widehat{D}^*P_C(x,y)(y^*) = \left\{ x^* \in X^*: \limsup_{\substack{(u,v) \to (x,y) \\ v \in P_C(u)}} \frac{\langle x^*, u-x \rangle - \langle y^*, v-y \rangle}{\|u-x\| + \|v-y\|} \leq 0 \right\}. \tag{2.10}$$

Notice that, if the considered Banach space is uniformly convex and uniformly smooth Banach spaces, then the metric projection operator $P_C$ becomes a single-valued continuous mapping. By (2.20) in [16], for any $x \in X$, in this case, we calculate the solutions of the Mordukhovich derivative of $P_C$ at $(x, P_C(x))$ as follows. For any $y^* \in X^*$, we have

$$\widehat{D}^*P_C(x, P_C(x))(y^*) := \widehat{D}^*P_C(x)(y^*)$$

$$= \left\{ x^* \in X^*: \limsup_{u \to x} \frac{\langle x^*, u-x \rangle - \langle y^*, P_C(u) - P_C(x) \rangle}{\|u-x\| + \|P_C(u) - P_C(x)\|} \leq 0 \right\}.$$

## 3. The real Banach space $l_1$

### 3.1. The metric projection operator in the real Banach space $l_1$

In this subsection, we consider the real Banach space $(l_1, \|\cdot\|_1)$ with dual space $(l_\infty, \|\cdot\|_\infty)$. It is well-known that $l_1$ is not uniformly convex and uniformly smooth. We will show that the metric projection $P_{r\mathbb{B}}$ is indeed a set-valued mapping. It causes that both the Gâteaux directional differentiability and Fréchet differentiability of $P_{r\mathbb{B}}$ are not defined in $l_1$. In this pace $l_1$, we can only study the Mordukhovich derivative of $P_{r\mathbb{B}}$ as a set-valued mapping.

We define some nonempty closed and convex cones in $l_1$ and $l_\infty$, respectively.

$$K_1 = \{x = (t_1, t_2, \ldots) \in l_1: t_n \geq 0, \text{ for all } n\},$$

$$K_1^+ = \{x = (t_1, t_2, \ldots) \in K_1 : t_n > 0, \text{ for all } n\},$$

and
$$K_\infty = \{\varphi = (u_1, u_2, \ldots) \in l_\infty : u_n \geq 0, \text{ for all } n\},$$

$$K_\infty^+ = \{\varphi = (u_1, u_2, \ldots) \in K_\infty : u_n > 0, \text{ for all } n\}.$$

$K_1$ and $K_\infty$ are the positive cones in $l_1$ and $l_\infty$, respectively. Let $\leqslant_1$ and $\leqslant_\infty$ be the two partial orders on $l_1$ and $l_\infty$ induced by the pointed closed and convex cones $K_1$ and $K_\infty$, respectively. We define the ordered intervals in $l_1$ and $l_\infty$. For any $x, y \in l_1$ with $x \leqslant_1 y$, and for any $u, v \in l_\infty$ with $u \leqslant_\infty v$, we write

$$[x, y]_{\leqslant_1} = \{z \in l_1 : x \leqslant_1 z \leqslant_1 y\},$$

$$[u, v]_{\leqslant_\infty} = \{w \in l_\infty : u \leqslant_\infty w \leqslant_\infty v\}.$$

For any $a > 0$, we define some nonempty convex subsets in $l_1$.

$$\Delta_a = \{x = (t_1, t_2, \ldots) \in K_1 : \|x\|_1 = \sum_{n=1}^\infty t_n = a\},$$

$$\Delta_a^+ = \{x = (t_1, t_2, \ldots) \in \Delta_a : t_n > 0, \text{ for all } n\}.$$

$\Delta_a$ is a nonempty closed and convex subset of $l_1$. In particular, $\Delta_1$ is simply denoted by $\Delta$ that is the simplex in $l_1$. For any real number $d$, we write $\beta_d = (d, d, \ldots) \in l_\infty$. In particular, we write $\beta = \beta_1$. If a point in $l_1$ or in $l_\infty$ has all entries 0 except one entry, then it is called an axis point in $l_1$ or in $l_\infty$, respectively. More precisely speaking, for a given point in $l_1$ or in $l_\infty$, if there is a positive number $m$ such that its $m^{th}$ entry equals $b \neq 0$ and all other entries are zero, then, it is denoted by $s(m, b)$.

As stated in section 2, the normalized duality mapping $J^*$ on $l_\infty$ is a set-valued mapping defined on $l_\infty$. We define a mapping $k^* : l_\infty \to l_1$, for $\varphi = (u_1, u_2, u_3, \ldots) \in l_\infty$, by

$$k^*(\varphi) = (\frac{u_1}{2}, \frac{u_2}{2^2}, \frac{u_3}{2^3}, \ldots) \in l_1.$$

$k^*$ has the following properties:

(a) $\|k^*(\varphi)\|_1 \leq \|\varphi\|_\infty$, for any $\varphi \in l_\infty$;
(b) $|\langle \varphi, k^*(\varphi) \rangle| \leq \|\varphi\|_\infty^2$, for any $\varphi \in l_\infty$;
(c) In particular, one has

$$k^*(\beta_a) \in J^*(\beta_a), \text{ for any } d \geq 0.$$

In the following proposition, we demonstrate that there are some points in $l_1$ that have non-singleton generalized identical sets. The proof of the following proposition is provided in the appendix (Partial of the results have been proved in [11]).

**Proposition 3.1.** For any $a > 0$, we have

(a) $J^{-1}(\beta_a) = \Delta_a$;

(b) $J(y) = \beta_a$, for any $y \in \Delta_a^+$;
(c) $J(y) \supsetneq \{\beta_a\}$, for any $y \in \Delta_a \setminus \Delta_a^+$;
(d) $\beta_a = \cap_{x \in \Delta_a} J(x)$;
(e) $\Im(y) = \Delta_a$, for every $y \in \Delta_a^+$;
(f) $\Im(y) \supsetneq \Delta_a$, for every $y \in \Delta_a \setminus \Delta_a^+$.

**Lemma 3.2**. *Let $r > 0$. For any $x \in l_1$ with $\|x\|_1 > r$, the set-valued projection $P_{r\mathbb{B}}$ satisfies*

(a) $\frac{r}{\|x\|_1} x \in P_{r\mathbb{B}}(x)$;
(b) *For any $j(x) \in J(x)$, we have*

$$\langle j(x), \frac{r}{\|x\|_1} x - y \rangle \geq 0, \text{ for all } y \in r\mathbb{B}.$$

*Proof.* Proof of (a). For any $x \in l_1$ with $\|x\|_1 > r$, it is clear that $\left\|\frac{r}{\|x\|_1} x\right\|_1 = r$. For any $z \in r\mathbb{B}$, we have

$$\|x - z\|_1 \geq \|x\|_1 - \|z\|_1 \geq \|x\|_1 - r.$$

On the other hand, we have

$$\left\|x - \frac{r}{\|x\|_1} x\right\|_1 = \|x\|_1 - r.$$

This implies $\frac{r}{\|x\|_1} x \in P_{r\mathbb{B}}(x)$. To prove part (b), we calculate

$$\langle j(x), \frac{r}{\|x\|_1} x - y \rangle$$

$$= \frac{r}{\|x\|_1} \|x\|_1^2 - \langle j(x), y \rangle$$

$$\geq \frac{r}{\|x\|_1} \|x\|_1^2 - \|j(x)\|_\infty \|y\|_1$$

$$= r\|x\|_1 - \|x\|_1 \|y\|_1$$

$$\geq 0, \text{ for all } y \in r\mathbb{B}. \qquad \square$$

**Proposition 3.3**. *Let $r > 0$. Let $x \in K_1$ with $\|x\|_1 > r$. we have*

(a) $P_{r\mathbb{B}}(x)$ *is a nonempty closed and convex subset of $\Delta_r$ satisfying*

$$P_{r\mathbb{B}}(x) = [\theta, x]_{\leqslant_1} \cap \Delta_r; \qquad (3.1)$$

(b) *For any $z \in K_1$ with $z \leqslant_1 x$, then*

$$P_{r\mathbb{B}}(z) \subseteq P_{r\mathbb{B}}(x);$$

(c) $P_{rB}(x)$ is a singleton if and only if $x$ is an axis point in $l_1$ satisfying

$$P_{rB}(s(m, \|x\|_1)) = s(m, r), \text{ for some positive integer } m. \tag{3.2}$$

*Proof.* We only prove (a). For this given $x \in K_1$ with $\|x\|_1 > r$ and for any $z \in rB$, we have

$$\|x - z\|_1 \geq \|x\|_1 - \|z\|_1 \geq \|x\|_1 - r. \tag{3.3}$$

Let $x = (t_1, t_2, \ldots) x \in K_1$. For any $y = (s_1, s_2, \ldots) \in [\theta, x]_{\leq_1} \cap \Delta_r$, we have

$$\|x - y\|_1 = \sum_{n=1}^{\infty} |t_n - s_n| = \sum_{n=1}^{\infty} (t_n - s_n) = \sum_{n=1}^{\infty} t_n - \sum_{n=1}^{\infty} s_n = \|x\|_1 - r. \tag{3.4}$$

By (3.3) and (3.4), it follows

$$y \in P_{rB}(x), \text{ for any } y \in [\theta, x]_{\leq_1} \cap \Delta_r.$$

This implies

$$[\theta, x]_{\leq_1} \cap \Delta_r \subseteq P_{rB}(x). \tag{3.5}$$

By Lemma 3.2, we have $\frac{r}{\|x\|_1} x \in [\theta, x]_{\leq_1} \cap \Delta_r \subseteq P_{rB}(x)$, which implies that $P_{rB}(x) \neq \emptyset$.

On the other hand of (3.5), we consider the following two cases.

Case 1. $y = (s_1, s_2, \ldots) \in [\theta, x]_{\leq_1} \cap rB$ and $y \notin \Delta_r$. That is $\|y\|_1 < r$. In this case, similar to the proof of (3.2), we have $\|x - y\|_1 > \|x\|_1 - r$, which implies

$$y \in [\theta, x]_{\leq_1} \cap rB \text{ and } y \notin \Delta_r \Longrightarrow y \notin P_{rB}(x). \tag{3.6}$$

Case 2. $y = (s_1, s_2, \ldots) \in rB$ and $y \notin [\theta, x]_{\leq_1}$. That is $\|y\|_1 \leq r$ and there is a positive number $m$ such that $s_m > t_m \geq 0$. In this case, similar to the proof of (3.2), we have

$$\|x - y\|_1$$
$$= \sum_{t_n \geq s_n} (t_n - s_n) + \sum_{t_n < s_n} (s_n - t_n)$$
$$= \sum_{t_n \geq s_n} t_n - \sum_{t_n \geq s_n} s_n + \sum_{t_n < s_n} s_n - \sum_{t_n < s_n} t_n$$
$$= \sum_{t_n \geq s_n} t_n - \sum_{t_n \geq s_n} s_n - \sum_{t_n < s_n} s_n + 2\sum_{t_n < s_n} s_n + \sum_{t_n < s_n} t_n - 2\sum_{t_n < s_n} t_n$$
$$= \sum_{n=1}^{\infty} t_n - \sum_{n=1}^{\infty} s_n + 2\sum_{t_n < s_n} (s_n - t_n)$$
$$> \|x\|_1 - \sum_{n=1}^{\infty} s_n$$
$$\geq \|x\|_1 - \sum_{n=1}^{\infty} |s_n|$$
$$\geq \|x\|_1 - r.$$

By (3.3) and (3.4), this implies $y \notin P_{rB}(x)$, that is,

$$y \in r\mathbb{B} \text{ and } y \notin [\theta, x]_{\leqslant_1} \implies y \notin P_{r\mathbb{B}}(x). \tag{3.5}$$

Then, (3.1) in this proposition is proved by (3.5), (3.6) and (3.7). Parts (b) follows from part (a) immediately

Proof of (c). For any positive integer $m$, we can check

$$[\theta, s(m, \|x\|_1)]_{\leqslant_1} \cap \Delta_r = \{s(m, r)\}.$$

By (3.1) in part (a), then (3.2) follows from this immediately.

Next, we prove that if $y$ is not an axis point, then $P_{r\mathbb{B}}(y)$ is a non-singleton by construction. Suppose that this given point $x = (t_1, t_2, \ldots)$ is not an axis point. Since $x \in K_1$ with $\|x\|_1 > r$, there are two positive integers $m$ and $k$ such that $t_k \geq t_m > 0$. By $0 < \frac{r}{\|x\|_1} < 1$, we take any $c$ with $1 > c > 0$ satisfying $(1+c)\frac{r}{\|x\|_1} < 1$. Then, we define $y_\lambda = (s_1, s_2, \ldots) \in \Delta_r$ as follows, for all $n$,

$$s_n = \begin{cases} \lambda, & \text{if } n = m, \\ \frac{r}{\|x\|_1}(t_m + t_k) - \lambda, & \text{if } n = k, \\ \frac{r}{\|x\|_1} t_n, & \text{if } n \neq m, k, \end{cases} \text{ for any } (1-c)\frac{r}{\|x\|_1} t_m \leq \lambda \leq \frac{r}{\|x\|_1} t_m.$$

Since $0 < \frac{r}{\|x\|_1} < 1$, one can check that

$$y_\lambda \in [\theta, x]_{\leqslant_1} \cap \Delta_r, \text{ for any } (1-c)\frac{r}{\|x\|_1} t_m \leq \lambda \leq \frac{r}{\|x\|_1} t_m.$$

By (3.1), $y_\lambda \in P_{r\mathbb{B}}(x)$. By $\frac{r}{\|x\|_1} x \in P_{r\mathbb{B}}(x)$, we have

$$y_\lambda \neq \frac{r}{\|x\|_1} x, \text{ for any } (1-c)\frac{r}{\|x\|_1} t_m \leq \lambda < \frac{r}{\|x\|_1} t_m.$$

This implies that $P_{r\mathbb{B}}(x)$ is not a singleton, which contains infinite many elements. □

### 3.2. The Mordukhovich derivatives of the metric projection operator in $l_1$

**Theorem 3.4.** *For any $r > 0$, the Mordukhovich derivatives of the set-valued metric projection $P_{r\mathbb{B}}: l_1 \rightrightarrows r\mathbb{B}$ has the following properties.*

(i) *For every $x \in r\mathbb{B}^\circ$ with $\|x\|_1 < r$, we have*

$$\widehat{D}^* P_{r\mathbb{B}}(x, x)(\varphi) = \{\varphi\}, \text{ for every } \varphi \in l_\infty.$$

(ii) *For every $x \in K_1^+$ with $\|x\|_1 > r$, we have*

(a) $\widehat{D}^* P_{r\mathbb{B}}\left(x, \frac{r}{\|x\|_1} x\right)(\theta^*) = \{\theta^*\};$

(b) $\widehat{D}^* P_{r\mathbb{B}}\left(x, \frac{r}{\|x\|_1} x\right)(\beta_d) = \emptyset, \text{ for any } d > 0;$

(c) $\widehat{D}^* P_{r\mathbb{B}}\left(x, \frac{r}{\|x\|_1} x\right)(J(x)) = \emptyset$.

*Proof.* Proof of (i). When the metric projection operator $P_{r\mathbb{B}}$ is restricted on the open subset $r\mathbb{B}^o$, it becomes the identity mapping on $r\mathbb{B}^o$ satisfying

$$\widehat{D}^* P_{r\mathbb{B}}(x, P_{r\mathbb{B}}(x))(\varphi) = \widehat{D}^* P_{r\mathbb{B}}(x, x)(\varphi), \text{ for every } \varphi \in l_\infty.$$

Then, the proof is very similar to the proof of part (i) of Theorem 3.1 in [16] and is omitted here.

Proof of (a) in (ii). Let $x = (t_1, t_2, \ldots) \in K_1^+$ with $\|x\|_1 := a > r$, It is clear that

$$\theta^* \in \widehat{D}^* P_{r\mathbb{B}}\left(x, \frac{r}{\|x\|_1} x\right)(\theta^*). \tag{3.6}$$

Let $\psi = (v_1, v_2, \ldots) \in l_\infty$. Suppose $\psi \neq \theta^*$. We have

$$\psi \in \widehat{D}^* P_{r\mathbb{B}}\left(x, \frac{r}{a} x\right)(\theta^*) \iff \limsup_{\substack{(u,v) \to (x, \frac{r}{a}x) \\ v \in P_{r\mathbb{B}}(u)}} \frac{\langle \psi, u-x \rangle - \langle \theta^*, v - \frac{r}{a}x \rangle}{\|u-x\|_1 + \|v - \frac{r}{a}x\|_1} \leq 0$$

$$\iff \limsup_{\substack{(u,v) \to (x, \frac{r}{a}x) \\ v \in P_{r\mathbb{B}}(u)}} \frac{\langle \psi, u-x \rangle}{\|u-x\|_1 + \|v - \frac{r}{a}x\|_1} \leq 0. \tag{3.7}$$

Notice that, for $x \in K_1^+$ with $\|x\|_1 > r$, there is $p$ with $0 < p < a - r$ satisfying

$$\|u\|_1 > r, \text{ for any } u \in K_1 \text{ with } \|u - x\|_1 < p.$$

By Lemma 3.2, this implies

$$\frac{r}{\|u\|_1} u \in P_{r\mathbb{B}}(u) \cap \Delta_r, \text{ for any } u \in K_1 \text{ with } \|u - x\|_1 < p. \tag{3.8}$$

It follows that, for any $u \in K_1$,

$$u \to x \iff \frac{r}{\|u\|_1} u \to \frac{r}{\|x\|_1} x.$$

This implies

$$\limsup_{\substack{(u,v) \to (x, \frac{r}{a}x) \\ v \in P_{r\mathbb{B}}(u)}} \frac{\langle \psi, u-x \rangle}{\|u-x\|_1 + \|v - \frac{r}{a}x\|_1} \geq \limsup_{\substack{(u, \frac{r}{\|u\|_1}u) \to (x, \frac{r}{a}x) \\ u \in K_1}} \frac{\langle \psi, u-x \rangle}{\|u-x\|_1 + \|\frac{r}{\|u\|_1}u - \frac{r}{a}x\|_1}. \tag{3.9}$$

By assumption that $\psi \neq \theta^*$. Define $\psi^+ \in l_\infty$ by

$$(\psi^+)_n = \begin{cases} v_n, & \text{if } v_n > 0, \\ 0, & \text{if } v_n \leq 0. \end{cases}$$

Define $\psi^- \in l_\infty$ by

$$(\psi^-)_n = \begin{cases} v_n, & \text{if } v_n < 0, \\ 0, & \text{if } v_n \geq 0. \end{cases}$$

Then, $\psi = \psi^+ + \psi^-$. By $\psi \neq \theta^*$, we have $\psi^+ \neq \theta^*$ or $\psi^- \neq \theta^*$, or both not $\theta^*$. Then, rest of the proof of (a) in (ii) is divided into two cases.

Case 1. $\psi^+ \neq \theta^*$. Then $k^*(\psi^+) \in K_1$ with $k^*(\psi^+) \neq \theta$. In limit (3.9), we take a line segment direction as $u = x + tk^*(\psi^+)$ for $t \downarrow 0$ with $0 < t < \frac{p}{\|k^*(\psi^+)\|_1}$. By $x \in K_1$, we have

$$x + tk^*(\psi^+) \in K_1, \text{ for all } 0 < t < \frac{p}{\|k^*(\psi^+)\|_1}.$$

By $k^*(\psi^+) \in K_1$ and $x \in K_1^+$, we have

$$|\|x + tk^*(\psi^+)\|_1 - \|x\|_1| = t\|k^*(\psi^+)\|_1, \text{ for all } 0 < t < \frac{p}{\|k^*(\psi^+)\|_1}.$$

By (3.9), and noticing $\langle \psi^+, k^*(\psi^+) \rangle > 0$, we have

$$\limsup_{\substack{(u, \frac{r}{\|u\|_1}u) \to (x, \frac{r}{a}x) \\ u \in K_1}} \frac{\langle \psi, u-x \rangle}{\|u-x\|_1 + \left\|\frac{r}{\|u\|_1}u - \frac{r}{a}x\right\|_1}$$

$$\geq \limsup_{\substack{t \downarrow 0 \\ t < \frac{p}{\|k^*(\psi^+)\|_1}}} \frac{\langle \psi, x + tk^*(\psi^+) - x \rangle}{\|x + tk^*(\psi^+) - x\|_1 + \left\|\frac{r}{\|x + tk^*(\psi^+)\|_1}(x + tk^*(\psi^+)) - \frac{r}{a}x\right\|_1}$$

$$= \limsup_{\substack{t \downarrow 0 \\ t < \frac{p}{\|k^*(\psi^+)\|_1}}} \frac{t\langle \psi, k^*(\psi^+) \rangle}{t\|k^*(\psi^+)\|_1 + r\left\|\left(\frac{1}{\|x + tk^*(\psi^+)\|_1} - \frac{1}{a}\right)x + \frac{trk^*(\psi^+)}{\|x + tk^*(\psi^+)\|_1}\right\|_1}$$

$$\geq \limsup_{\substack{t \downarrow 0 \\ t < \frac{p}{\|k^*(\psi^+)\|_1}}} \frac{t\langle \psi^+, k^*(\psi^+) \rangle}{t\|k^*(\psi^+)\|_1 + r\|x\|_1\left|\frac{1}{\|x + tk^*(\psi^+)\|_1} - \frac{1}{a}\right| + \frac{tr\|k^*(\psi^+)\|_1}{\|x + tk^*(\psi^+)\|_1}}$$

$$= \limsup_{\substack{t \downarrow 0 \\ t < \frac{p}{\|k^*(\psi^+)\|_1}}} \frac{t\langle \psi^+, k^*(\psi^+) \rangle}{t\|k^*(\psi^+)\|_1 + r\|x\|_1 \frac{1}{\|x + tk^*(\psi^+)\|_1 \|x\|_1}|\|x + tk^*(\psi^+)\|_1 - \|x\|_1| + \frac{tr\|k^*(\psi^+)\|_1}{\|x + tk^*(\psi^+)\|_1}}$$

$$= \limsup_{\substack{t \downarrow 0 \\ t < \frac{p}{\|k^*(\psi^+)\|_1}}} \frac{t\langle \psi^+, k^*(\psi^+) \rangle}{t\|k^*(\psi^+)\|_1 + \frac{r}{\|x + tk^*(\psi^+)\|_1}|\|x + tk^*(\psi^+)\|_1 - \|x\|_1| + \frac{tr\|k^*(\psi^+)\|_1}{\|x + tk^*(\psi^+)\|_1}}$$

$$= \limsup_{\substack{t \downarrow 0 \\ t < \frac{p}{\|k^*(\psi^+)\|_1}}} \frac{t\langle \psi^+, k^*(\psi^+) \rangle}{t\|k^*(\psi^+)\|_1 + \frac{rt\|k^*(\psi^+)\|_1}{\|x + tk^*(\psi^+)\|_1} + \frac{tr\|k^*(\psi^+)\|_1}{\|x + tk^*(\psi^+)\|_1}}$$

$$= \limsup_{\substack{t \downarrow 0 \\ t < \frac{p}{\|k^*(\psi^+)\|_1}}} \frac{\langle \psi^+, k^*(\psi^+) \rangle}{\|k^*(\psi^+)\|_1 + \frac{r\|k^*(\psi^+)\|_1}{\|x+tk^*(\psi^+)\|_1} + \frac{r\|k^*(\psi^+)\|_1}{\|x+tk^*(\psi^+)\|_1}}$$

$$= \frac{\langle \psi^+, k^*(\psi^+) \rangle}{\|k^*(\psi^+)\|_1 + \frac{2r\|k^*(\psi^+)\|_1}{\|x\|_1}}$$

$$> 0.$$

This shows that

$$\psi^+ \neq \theta^* \Longrightarrow \psi \notin \widehat{D}^* P_{r\mathbb{B}}\left(x, \frac{r}{a}x\right)(\theta^*). \tag{3.10}$$

Case 2. $\psi^- \neq \theta^*$. In this case, there is a positive integer $m$ such that $\{n \leq m : v_n < 0\} \neq \emptyset$. Define $\psi_m^- \in l_1 \cap l_\infty$, for all $n$ by

$$(\psi_m^-)_n = \begin{cases} v_n, & \text{if } n \leq m \text{ and } v_n < 0, \\ 0, & \text{otherwise.} \end{cases}$$

Recall that $p$ is taken to be positive with $0 < p < a - r$ satisfying

$$\|u\|_1 > r, \text{ for any } u \in K_1 \text{ with } \|u - x\|_1 < p.$$

Since $x \in K_1^+$, then $\min\{t_n : n \leq m\} > 0$. Let $b = \min\{\min\{t_n : n \leq m\}, p\}$, in which $p$ is defined in (3.8). In this case, we have

$$x + t\psi_m^- \in K_1 \text{ and } \|x + t\psi_m^-\|_1 > r, \text{ for all } 0 < t < \frac{b}{m\|\psi\|_\infty}.$$

Then, in the limit (3.9), we take a line segment direction as $u_t = x + t\psi_m^-$, for $t \downarrow 0$ with $0 < t < \frac{b}{m\|\psi\|_\infty}$. We have

$$\langle \psi, \psi_m^- \rangle = \langle \psi_m^-, \psi_m^- \rangle > 0,$$

and

$$|\|x + t\psi_m^-\|_1 - \|x\|_1| = t\|\psi_m^-\|_1, \text{ for } 0 < t < \frac{b}{m\|\psi\|_\infty}.$$

By (3.9), this implies

$$\limsup_{\substack{\left(u, \frac{r}{\|u\|_1}u\right) \to \left(x, \frac{r}{a}x\right) \\ u \in K_1}} \frac{\langle \psi, u - x \rangle}{\|u - x\|_1 + \left\|\frac{r}{\|u\|_1}u - \frac{r}{a}x\right\|_1}$$

$$\geq \limsup_{\left(u_t, \frac{r}{\|u_t\|_1}u_t\right) \to \left(x, \frac{r}{a}x\right)} \frac{\langle \psi, u_t - x \rangle}{\|u_t - x\|_1 + \left\|\frac{r}{\|u_t\|_1}u_t - \frac{r}{a}x\right\|_1}$$

$$\begin{aligned}
&= \limsup_{\substack{t\downarrow 0 \\ t<\frac{b}{m\|\psi\|_\infty}}} \frac{\langle \psi,\ x+t\psi_m^- - x\rangle}{\|x+t\psi_m^- - x\|_1 + \left\|\frac{r}{\|x+t\psi_m^-\|_1}(x+t\psi_m^-) - \frac{r}{a}x\right\|_1} \\
&= \limsup_{\substack{t\downarrow 0 \\ t<\frac{b}{m\|\psi\|_\infty}}} \frac{t\langle \psi_m^-,\ \psi_m^-\rangle}{\|x+t\psi_m^- - x\|_1 + \left\|\frac{r}{\|x+t\psi_m^-\|_1}(x+t\psi_m^-) - \frac{r}{a}x\right\|_1} \\
&\geq \limsup_{\substack{t\downarrow 0 \\ t<\frac{b}{m\|\psi\|_\infty}}} \frac{t\langle \psi_m^-,\ \psi_m^-\rangle}{t\|\psi_m^-\|_1 + \frac{r}{\|x+t\psi_m^-\|_1}\big|\|x+t\psi_m^-\|_1 - \|x\|_1\big| + \frac{tr\|\psi_m^-\|_1}{\|x+\psi_m^-\|_1}} \\
&= \limsup_{\substack{t\downarrow 0 \\ t<\frac{b}{m\|\psi\|_\infty}}} \frac{t\langle \psi_m^-,\ \psi_m^-\rangle}{t\|\psi_m^-\|_1 + \frac{rt\|\psi_m^-\|_1}{\|x+t\psi_m^-\|_1} + \frac{tr\|\psi_m^-\|_1}{\|x+t\psi_m^-\|_1}} \\
&= \limsup_{\substack{t\downarrow 0 \\ t<\frac{b}{m\|\psi\|_\infty}}} \frac{\langle \psi_m^-,\ \psi_m^-\rangle}{\|\psi_m^-\|_1 + \frac{r\|\psi_m^-\|_1}{\|x+t\psi_m^-\|_1} + \frac{r\|\psi_m^-\|_1}{\|x+t\psi_m^-\|_1}} \\
&= \frac{\langle \psi_m^-,\ \psi_m^-\rangle}{\|\psi_m^-\|_1 + \frac{2r\|\psi_m^-\|_1}{\|x\|_1}} \\
&> 0.
\end{aligned}$$

This shows that

$$\psi^- \neq \theta^* \Longrightarrow \psi \notin \widehat{D}^* P_{r\mathbb{B}}\left(x, \frac{r}{a}x\right)(\theta^*). \tag{3.11}$$

By (3.10) and (3.11), it follows that

$$\psi \neq \theta^* \Longrightarrow \psi \notin \widehat{D}^* P_{r\mathbb{B}}\left(x, \frac{r}{a}x\right)(\theta^*). \tag{3.12}$$

By (3.6) and (3.12), we obtain

$$\widehat{D}^* P_{r\mathbb{B}}\left(x, \frac{r}{a}x\right)(\theta^*) = \{\theta^*\}.$$

Proof of (b) in (ii). For $\frac{r}{\|x\|_1}x \in P_{r\mathbb{B}}(x)$, by Proposition 3.2, we have

$$\frac{r}{\|x\|_1}x \in [\theta, x]_{\leqslant_1} \cap \Delta_r.$$

Let $d > 0$. One has

$$\langle \beta_d, \frac{r}{\|x\|_1}x\rangle = dr. \tag{3.13}$$

For any $\psi \in l_\infty$, we have

$$\psi \in \widehat{D}^*P_{r\mathbb{B}}\left(x, \frac{r}{\|x\|_1}x\right)(\beta_d) \Leftrightarrow \limsup_{\substack{(u,v)\to\left(x,\frac{r}{\|x\|_1}x\right) \\ v\in P_{r\mathbb{B}}(u)}} \frac{\langle\psi,u-x\rangle - \langle\beta_d, v-\frac{r}{\|x\|_1}x\rangle}{\|u-x\|_1 + \|v-\frac{r}{\|x\|_1}x\|_1} \leq 0. \quad (3.14)$$

At first, we prove

$$\theta^* \notin \widehat{D}^*P_{r\mathbb{B}}\left(x, \frac{r}{\|x\|_1}x\right)(\beta_d). \quad (3.15)$$

Assume, be the way of contradiction, that $\theta^* \in \widehat{D}^*P_{r\mathbb{B}}\left(x, \frac{r}{\|x\|_1}x\right)(\beta_d)$. By (3.14), we have

$$\theta^* \in \widehat{D}^*P_{r\mathbb{B}}\left(x, \frac{r}{\|x\|_1}x\right)(\beta_d) \Leftrightarrow \limsup_{\substack{(u,v)\to\left(x,\frac{r}{\|x\|_1}x\right) \\ v\in P_{r\mathbb{B}}(u)}} \frac{-\langle\beta_d, v-\frac{r}{\|x\|_1}x\rangle}{\|u-x\|_1 + \|v-\frac{r}{\|x\|_1}x\|_1} \leq 0$$

$$\Leftrightarrow \limsup_{\substack{(u,v)\to\left(x,\frac{r}{\|x\|_1}x\right) \\ v\in P_{r\mathbb{B}}(u)}} \frac{\langle\beta_d, \frac{r}{\|x\|_1}x-v\rangle}{\|u-x\|_1 + \|v-\frac{r}{\|x\|_1}x\|_1} \leq 0. \quad (3.16)$$

Notice that, $x = (t_1, t_2, \ldots) \in K_1^+$. We have $t_n > 0$, for all $n$ satisfying $t_n \to 0$, as $n \to \infty$. Then, in the limit (3.16), we take a sequence direction as $u_n = x + s(n, -\left(t_n + \frac{1}{n}\right))$, for $n \to \infty$. It is clear that $u_n = x + s(n, -\left(t_n + \frac{1}{n}\right)) \notin K_1$, for all $n$ satisfying $u_n \to x$, as $n \to \infty$. By Lemma 3.1, we have $\frac{r}{\|u_n\|_1}u_n \in P_{r\mathbb{B}}(u_n)$. We calculate

$$\|u_n\|_1 = \left\|x + s\left(n, -\left(t_n + \frac{1}{n}\right)\right)\right\|_1 = \|x\|_1 - t_n + \frac{1}{n} \to \|x\|_1, \text{ as } n \to \infty.$$

By $u_n \to x$, as $n \to \infty$, this implies

$$\frac{r}{\|u_n\|_1}u_n \to \frac{r}{\|x\|_1}x, \text{ as } n \to \infty. \quad (3.17)$$

And

$$\langle\beta_d, \frac{r}{\|u_n\|_1}u_n\rangle = \frac{r}{\|u_n\|_1}\langle\beta_d, u_n\rangle = \frac{dr}{\|x\|_1 - t_n + \frac{1}{n}}\left(\|x\|_1 - t_n - \frac{1}{n}\right), \text{ for all } n. \quad (3.18)$$

As taking the limit for $n$, we can assume $\|x\|_1 - t_n - 1/n > 0$, for $n \gg 1$. In the limit in (3.16), by (3.13), (3.17) and (3.18), we have

$$\limsup_{\substack{(u,v)\to\left(x,\frac{r}{\|x\|_1}x\right) \\ v\in P_{r\mathbb{B}}(u)}} \frac{\langle\beta_d, \frac{r}{\|x\|_1}x-v\rangle}{\|u-x\|_1 + \|v-\frac{r}{\|x\|_1}x\|_1}$$

$$
\begin{aligned}
&\geq \limsup_{\left(u_n, \frac{r}{\|u_n\|_1} u_n\right) \to \left(x, \frac{r}{\|x\|_1} x\right)} \frac{\langle \beta_d, \frac{r}{\|x\|_1} x - \frac{r}{\|u_n\|_1} u_n \rangle}{\|u_n - x\|_1 + \left\| \frac{r}{\|u_n\|_1} u_n - \frac{r}{\|x\|_1} x \right\|_1} \\
&= \limsup_{\left(u_n, \frac{r}{\|u_n\|_1} u_n\right) \to \left(x, \frac{r}{\|x\|_1} x\right)} \frac{\langle \beta_d, \frac{r}{\|x\|_1} x \rangle - \frac{r}{\|u_n\|_1} \langle \beta_b, u_n \rangle}{\|u_n - x\|_1 + \left\| \frac{r}{\|u_n\|_1} u_n - \frac{r}{\|x\|_1} x \right\|_1} \\
&= \limsup_{\left(u_n, \frac{r}{\|u_n\|_1} u_n\right) \to \left(x, \frac{r}{\|x\|_1} x\right)} \frac{dr - \frac{br}{\|x\|_1 - t_n + \frac{1}{n}} \left( \|x\|_1 - t_n - \frac{1}{n} \right)}{\left( t_n + \frac{1}{n} \right) + r \left\| \frac{1}{\|u_n\|_1} u_n - \frac{1}{\|x\|_1} x \right\|_1} \\
&= \limsup_{\left(u_n, \frac{r}{\|u_n\|_1} u_n\right) \to \left(x, \frac{r}{\|x\|_1} x\right)} \frac{\frac{2dr}{\|x\|_1 - t_n + \frac{1}{n}} \frac{1}{n}}{\left( t_n + \frac{1}{n} \right) + r \left\| \frac{1}{\|u_n\|_1} u_n - \frac{1}{\|x\|_1} x \right\|_1} \\
&= \limsup_{\left(u_n, \frac{r}{\|u_n\|_1} u_n\right) \to \left(x, \frac{r}{\|x\|_1} x\right)} \frac{\frac{2dr}{\|x\|_1 - t_n + \frac{1}{n}} \frac{1}{n}}{\left( t_n + \frac{1}{n} \right) + r \left\| \frac{1}{\|x\|_1 - t_n + \frac{1}{n}} \left( x + s(n, -(t_n + \frac{1}{n})) \right) - \frac{1}{\|x\|_1} x \right\|_1} \\
&\geq \limsup_{\left(u_n, \frac{r}{\|u_n\|_1} u_n\right) \to \left(x, \frac{r}{\|x\|_1} x\right)} \frac{\frac{2dr}{\|x\|_1 - t_n + \frac{1}{n}} \frac{1}{n}}{\left( t_n + \frac{1}{n} \right) + r \left| \frac{1}{\|x\|_1 - t_n + \frac{1}{n}} - \frac{1}{\|x\|_1} \right| \|x\|_1 + r \left\| \frac{s(n, -(t_n + \frac{1}{n}))}{\|x\|_1 - t_n + \frac{1}{n}} \right\|_1} \\
&= \limsup_{\left(u_n, \frac{r}{\|u_n\|_1} u_n\right) \to \left(x, \frac{r}{\|x\|_1} x\right)} \frac{\frac{2dr}{\|x\|_1 - t_n + \frac{1}{n}} \frac{1}{n}}{\left( t_n + \frac{1}{n} \right) + \frac{r|t_n - \frac{1}{n}|}{\|x\|_1 - t_n + \frac{1}{n}} + \frac{r(t_n + \frac{1}{n})}{\|x\|_1 - t_n + \frac{1}{n}}} \\
&= \limsup_{\left(u_n, \frac{r}{\|u_n\|_1} u_n\right) \to \left(x, \frac{r}{\|x\|_1} x\right)} \frac{\frac{2dr}{\|x\|_1 - t_n + \frac{1}{n}}}{nt_n + \frac{r|nt_n - 1|}{\|x\|_1 - t_n + \frac{1}{n}} + \frac{r(nt_n + 1)}{\|x\|_1 - t_n + \frac{1}{n}}} \\
&= \limsup_{\left(u_n, \frac{r}{\|u_n\|_1} u_n\right) \to \left(x, \frac{r}{\|x\|_1} x\right)} \frac{\frac{2dr}{\|x\|_1}}{\frac{2r}{\|x\|_1}} \\
&= d.
\end{aligned}
$$

Where, we use the property that $x \in K_1^+$ implies $nt_n \to 0$, as $n \to \infty$. Since $c > 0$, by (3.16), this implies

$$\theta^* \notin \widehat{D}^* P_{r\mathbb{B}}\left(x, \frac{r}{\|x\|_1} x\right)(\beta_d). \tag{3.19}$$

For any $\psi = (\lambda_1, \lambda_2, \ldots) \in l_\infty$ with $\psi \neq \theta^*$. There is a positive number $m$ satisfying $\lambda_m \neq 0$. Then, in the limit (3.14), we take a segment direction as $u_t = x + s(n, t\lambda_m)$, for $t \downarrow 0$. Since $x = (t_1, t_2, \ldots) \in K_1^+$, it follows that

$$u_t = x + s(n, t\lambda_m) \in K_1^+ \text{ and } \frac{r}{\|u_t\|_1} u_t \in K_1^+ \text{ for all } 0 < t < \frac{t_m}{|\lambda_m|}.$$

One has that

$$\frac{r}{\|u_t\|_1} u_t \to \frac{r}{\|x\|_1} x, \text{ as } t \downarrow 0.$$

From (3.14), this implies

$$\limsup_{\substack{(u,v) \to \left(x, \frac{r}{\|x\|_1} x\right) \\ v \in P_{r\mathbb{B}}(u)}} \frac{\langle \psi, u-x \rangle - \langle \beta_d, v - \frac{r}{\|x\|_1} x \rangle}{\|u-x\|_1 + \|v - \frac{r}{\|x\|_1} x\|_1}$$

$$\geq \limsup_{\substack{\left(u_t, \frac{r}{\|u_t\|_1} u_t\right) \to \left(x, \frac{r}{\|x\|_1} x\right) \\ t < \frac{t_m}{|\lambda_m|}}} \frac{\langle \psi, u_t - x \rangle - \langle \beta_d, \frac{r}{\|u_t\|_1} u_t - \frac{r}{\|x\|_1} x \rangle}{\|u_t - x\|_1 + \left\|\frac{r}{\|u_t\|_1} u_t - \frac{r}{\|x\|_1} x\right\|_1}$$

$$= \limsup_{\substack{\left(u_t, \frac{r}{\|u_t\|_1} u_t\right) \to \left(x, \frac{r}{\|x\|_1} x\right) \\ t < \frac{t_m}{|\lambda_m|}}} \frac{\langle \psi, s(n, t\lambda_m) \rangle - (dr - dr)}{\|s(n, t\lambda_m)\|_1 + \left\|\frac{r}{\|x + s(n, t\lambda_m)\|_1}(x + s(n, t\lambda_m)) - \frac{r}{\|x\|_1} x\right\|_1}$$

$$= \limsup_{\substack{\left(u_t, \frac{r}{\|u_t\|_1} u_t\right) \to \left(x, \frac{r}{\|x\|_1} x\right) \\ t < \frac{t_m}{|\lambda_m|}}} \frac{t\lambda_m^2}{t|\lambda_m| + \left\|\frac{r}{\|x\|_1 + t\lambda_m}(x + s(n, t\lambda_m)) - \frac{r}{\|x\|_1} x\right\|_1}$$

$$\geq \limsup_{\substack{\left(u_t, \frac{r}{\|u_t\|_1} u_t\right) \to \left(x, \frac{r}{\|x\|_1} x\right) \\ t < \frac{t_m}{|\lambda_m|}}} \frac{t\lambda_m^2}{t|\lambda_m| + \frac{r}{\|x\|_1 + t\lambda_m} |\|x\|_1 + t\lambda_m - \|x\|_1| + r \left\|\frac{s(n, t\lambda_m)}{\|x\|_1 + t\lambda_m}\right\|_1}$$

$$= \limsup_{\substack{\left(u_t, \frac{r}{\|u_t\|_1} u_t\right) \to \left(x, \frac{r}{\|x\|_1} x\right) \\ t < \frac{t_m}{|\lambda_m|}}} \frac{t\lambda_m^2}{t|\lambda_m| + \frac{2rt|\lambda_m|}{\|x\|_1 + t\lambda_m}}$$

$$= \limsup_{\substack{\left(u_t, \frac{r}{\|u_t\|_1} u_t\right) \to \left(x, \frac{r}{\|x\|_1} x\right) \\ t < \frac{t_m}{|\lambda_m|}}} \frac{\lambda_m^2}{|\lambda_m| + \frac{2r|\lambda_m|}{\|x\|_1 + t\lambda_m}}$$

$$= \frac{|\lambda_m|}{1 + \frac{2r}{\|x\|_1}}$$

$$> 0.$$

By (3.14), this implies

$$\psi \neq \theta^* \implies \psi \notin \widehat{D}^* P_{r\mathbb{B}}\left(x, \frac{r}{\|x\|_1}x\right)(\beta_d). \tag{3.20}$$

By (3.19) and (3.20), we obtain

$$\widehat{D}^* P_{r\mathbb{B}}\left(x, \frac{r}{\|x\|_1}x\right)(\beta_b) = \emptyset, \text{ for any } d > 0.$$

Proof of (c) in (ii). For every $x \in K_1^+$, by part (b) of Proposition 3.1, $J(x)$ is a singleton satisfying $J(x) = \{\beta_{\|x\|_1}\}$. In addition, if $\|x\|_1 > r$, by part (b) in (ii) of this theorem, we have

$$\widehat{D}^* P_{r\mathbb{B}}\left(x, \frac{r}{\|x\|_1}x\right)(J(x)) = \widehat{D}^* P_{r\mathbb{B}}\left(x, \frac{r}{\|x\|_1}x\right)(\beta_{\|x\|_1}) = \emptyset. \qquad \square$$

**Remarks 3.5.** For any $r > 0$, in Theorem 3.4, we investigate some properties of Mordukhovich derivatives of the set-valued metric projection $P_{r\mathbb{B}}: l_1 \rightrightarrows r\mathbb{B}$, for $x \in K_1^+$ with $\|x\|_1 > r$. We can similarly study the case for $x \in -K_1^+$ with $\|x\|_1 > r$.

## 4. Banach space $c$

### 4.1. The set-valued metric projection from $c$ to $c_0$

Let $(c, \|\cdot\|)$ be the Banach space of all convergent sequences of real numbers (see [6]) with the supremum norm $\|\cdot\|$. For any $x = \{t_n\} \in c$, the norm $\|x\|$ of $x$ is defined by

$$\|x\| = \sup_{1 \leq n < \infty} |t_n|.$$

The dual space of $c$ is $c^* = l_1$. The norm of $c^* = l_1$ is denoted by $\|\cdot\|_1$ as used in the previous section. Since all elements in both $c$ and $c^*$ are sequences of real numbers, to make the distinction between the elements in $c$ and $c^*$, we use English letters to name the elements in $c$ and Greek letters to name the elements in $c^* = l_1$. Let $\theta$ and $\theta^*$ denote the origins in $c$ and $c^*$, respectively, which are identical with $\theta = \theta^* = (0, 0, \ldots)$. Let $\langle \cdot, \cdot \rangle$ denote the real canonical pairing between $c^*$ and $c$. For any $\varphi = (q_0, q_1, q_2, \ldots) \in c^* = l_1$, and $x = \{t_n\} \in c$, we have (see IV 2.7, II 4.31 and II 4.33 in [6]).

$$\langle \varphi, x \rangle = q_0 \lim_{n \to \infty} t_n + \sum_{n=1}^{\infty} q_n t_n.$$

Let $c_0$ be the closed subspace of $c$ that contains all convergent real sequences with their limit zero. In this section, we investigate the properties of Mordukhovich derivatives of the set-valued metric projection $P_{c_0}: c \rightrightarrows c_0$.

In the previous section, we defined $K_1$ and $K_1^+$ in $l_1$. Similar to those sets, we define the positive cones in $(c, \|\cdot\|)$ by

$$K_c = \{x = (t_1, t_2, \ldots) \in c: t_n \geq 0, \text{ for all } n\},$$

$$K_c^+ = \{x = (t_1, t_2, \ldots) \in K_c : t_n > 0, \text{ for all } n\}.$$

$K_c$ is the positive cone in $c$. Let $\leqslant_c$ be the partial order on $c$ induced by this pointed closed and convex cone $K_c$. We define the $\leqslant_c$-intervals in $c$. For any $x, y \in c$ with $x \leqslant_c y$, we write

$$[x, y]_{\leqslant_c} = \{z \in c : x \leqslant_c z \leqslant_c y\}.$$

For an arbitrary given real number $b$, let $\beta_b = (b, b, \ldots) \in c$, which is used in the previous section. In particular, $\beta_0 = \theta = \theta^*$. However, the $c$-norm of $\beta_d$ is $|d|$, which equals to its $l_\infty$-norm.

For any $r > 0$, let $\Delta_r$ be defined as in the previous section. Let $J$ be the normalized duality mapping from $c$ to $c^* = l_1$.

**Lemma 4.1.** *For any $r > 0$, let $\beta_r = (r, r, \ldots) \in c \backslash c_0$. Then*

$$J(\beta_r) = \Delta_r \subseteq l_1.$$

*Proof.* Since $\|\beta_r\| = r$, for any $\varphi = (q_0, q_1, q_2, \ldots) \in J(\beta_r) \subseteq l_1$, we have

$$r^2 = \|\beta_r\|^2 = \langle \varphi, \beta_r \rangle = \|\varphi\|_1^2.$$

It follows that

$$r^2 = \langle \varphi, \beta_r \rangle = q_0 \lim_{n \to \infty} r + \sum_{n=1}^{\infty} q_n r = r q_0 + r \sum_{n=1}^{\infty} q_n = r^2.$$

This implies,

$$\sum_{n=0}^{\infty} q_n = r.$$

For $\varphi = (q_0, q_1, q_2, \ldots) \in J(\beta_r) \subseteq l_1$, we have

$$\sum_{n=0}^{\infty} q_n = r = \|\beta_r\| = \|\varphi\|_1 = \sum_{n=0}^{\infty} |q_n|.$$

This implies that

$$0 \leq q_n \leq r, \text{ for } n = 0, 1, 2, \ldots.$$

Hence, $\varphi \in \Delta_r \subseteq l_1$. On the other hand, for any $\varphi = (q_0, q_1, q_2, \ldots) \in \Delta_r \subseteq l_1$, we calculate

$$\langle \varphi, \beta_r \rangle = q_0 \lim_{n \to \infty} r + \sum_{n=1}^{\infty} q_n r = r q_0 + r \sum_{n=1}^{\infty} q_n = r \sum_{n=0}^{\infty} q_n = r^2 = \|\varphi\|_1^2 = r^2 = \|\beta_r\|^2.$$

It follows that $\varphi \in J(\beta_r) \subseteq l_1$. Hence, we have

$$J(\beta_r) = \Delta_r \subseteq l_1, \text{ for } \beta_r = (r, r, r, \ldots) \in c \backslash c_0. \qquad \square$$

For any $x = (t_1, t_2, \ldots) \in c$, we write

$$L(x) = \lim_{n \to \infty} t_n \quad \text{and} \quad h(x) = x - \beta_{L(x)}.$$

**Lemma 4.2.** *For $x = (t_1, t_2, ...) \in c$, we have*

(a) *If $x \in c_0$, then $L(x) = 0$ and $h(x) = x$;*
(b) *$L(x+\beta_b) = L(x) + b$ and $h(x+\beta_b) = h(x)$, for any real number $b$.*

*Proof.* The proof of this lemma is straightforward and it is omitted here. □

Since the Banach space $c$ is not uniformly convex and uniformly smooth, the metric projection $P_{c_0}$ is not necessarily to be single-valued. As a matter of fact, the following proposition shows that the metric projection $P_{c_0} \colon c \rightrightarrows c_0$ is indeed a set-valued mapping, and it precisely describes its solutions.

**Proposition 4.3.** *For any $x = (t_1, t_2, ...) \in c$, we have*

$$P_{c_0}(x) = \{y \in c_0 \colon \|x - y\| = |L(x)|\}.$$

*In particular, we have*

(a) $h(x) := x - \beta_{L(x)} \in P_{c_0}(x)$;
(b) $P_{c_0}(\beta_b) = \{y \in c_0 \colon \|\beta_b - y\| = |b|\}$, *for any $b \neq 0$.*

*Proof.* If $x \in c_0$, then, $L(x) = 0$ and $P_{c_0}(x) = x$, which is a singleton. So, this proposition is satisfied for $x \in c_0$. Next, for an arbitrarily given $x = (t_1, t_2, ...) \in c \setminus c_0$, by definition, we have that $L(x) = \lim_{n \to \infty} t_n \neq 0$. Then, for any $z = (s_1, s_2, ...) \in c_0$, by $\lim_{n \to \infty} s_n = 0$, it follows that

$$\|x - z\| = \sup_{1 \leq n < \infty} |t_n - s_n| \geq |L(x)|, \text{ for any } z \in c_0. \tag{4.1}$$

(4.1) implies that

$$\{y \in c_0 \colon \|x - y\| = |L(x)|\} \subseteq P_{c_0}(x). \tag{4.2}$$

On the other hand, for any $y \in P_{c_0}(x)$. By definition, we have

$$\|x - y\| \leq \|x - z\|, \text{ for any } z \in c_0. \tag{4.3}$$

For the given $x = (t_1, t_2, ...) \in c \setminus c_0$ and $L(x) = \lim_{n \to \infty} t_n$ satisfying $L(x) \neq 0$. Take $h(x) = x - \beta_{L(x)} \in c_0$. Substituting $z$ by $h(x)$ in (4.3), we have

$$\|x - y\| \leq \|x - h(x)\| = \|\beta_{L(x)}\| = |L(x)|. \tag{4.4}$$

Since $y \in c_0$, by (4.1) we have $\|x - y\| \geq |L(x)|$. By (4.4), this implies $\|x - y\| = |L(x)|$. Then,

$$P_{c_0}(x) \subseteq \{y \in c_0 \colon \|x - y\| = |L(x)|\}. \tag{4.5}$$

By (4.2) and (4.5), this proposition is proved. □

It is clear that

$$u \to x \implies h(u) \to h(x), \text{ for } u \in c. \tag{4.6}$$

The following theorem provides the properties of Mordukhovich derivatives of the set-valued metric projection $P_{c_0}: c \rightrightarrows c_0$. In the following theorem, for $q_0 \in \mathbb{R}$, the axis point $s(0, q_0) = (q_0, 0, 0, \ldots)$, and plays a crucial role. Notice that, every axis point is in $c \cap c^*$.

### 4.2. The Mordukhovich derivatives of the set-valued metric projection in $c$

In this subsection, we use the properties of the set-valued metric projection $P_{c_0}: c \rightrightarrows c_0$ proved in the previous subsection to investigate the solutions of the Mordukhovich derivatives of $P_{c_0}$. We will find that the axis points of the form $s(0, q_0) = (q_0, 0, 0, \ldots)$ play crucial role for sovling the Mordukhovich derivatives of $P_{c_0}$. Since for every $x \in c$, we have $h(x) \in P_{c_0}(x)$. In the following theorem, we investigate the solutions of the Mordukhovich derivatives of $P_{c_0}$ at point $(x, h(x))$, for any $x \in c$.

**Theorem 4.4**. *For every $x \in c$, we have*

(i) *If $\varphi = s(0, q_0) \in c^*$, for some $q_0 \in \mathbb{R}$, then*

$$\widehat{D}^* P_{c_0}(x, h(x))(\varphi) = \{\theta^*\};$$

(ii) *For $\varphi \in c^*$, if $\varphi \neq s(0, b)$, for any $b \in \mathbb{R}$, then*

$$\widehat{D}^* P_{c_0}(x, h(x))(\varphi) = \emptyset.$$

*Proof.* Proof of (i). Let $x \in c$. For any $\varphi = s(0, q_0) \in c^*$, for some $q_0 \in \mathbb{R}$, we have

$$\theta^* \in \widehat{D}^* P_{c_0}(x, h(x))(\varphi) \iff \limsup_{\substack{(u,v) \to (x, h(x)) \\ v \in P_{c_0}(u)}} \frac{\langle \theta^*, u-x \rangle - \langle \varphi, v - h(x) \rangle}{\|u-x\| + \|v - h(x)\|} \leq 0. \tag{4.7}$$

For any $u \in c$ and for any $v \in P_{c_0}(u)$, we have $v \in c_0$. Since $h(x) \in c_0$, it is clear that

$$\langle \varphi, v - h(x) \rangle = \langle s(0, q_0), v - h(x) \rangle = q_0 L(v - h(x)) = 0. \tag{4.8}$$

This implies

$$\limsup_{\substack{(u,v) \to (x, h(x)) \\ v \in P_{c_0}(u)}} \frac{\langle \theta^*, u-x \rangle - \langle \varphi, v - h(x) \rangle}{\|u-x\| + \|v - h(x)\|}$$

$$= \limsup_{\substack{(u,v) \to (x, h(x)) \\ v \in P_{c_0}(u)}} \frac{-\langle s(0, q_0), v - h(x) \rangle}{\|u-x\| + \|v - h(x)\|}$$

$$= \limsup_{\substack{(u,v) \to (x, h(x)) \\ v \in P_{c_0}(u)}} \frac{0}{\|u-x\| + \|v - h(x)\|} = 0.$$

By (4.7), this implies

$$\theta^* \in \widehat{D}^* P_{c_0}(x, h(x))(\varphi), \text{ for } \varphi = s(0, q_0) \in c^* \text{ with } q_0 \in \mathbb{R}. \tag{4.9}$$

For any $\psi = (p_0, p_1, p_2, \dots) \in c^*$ with $\psi \neq \theta^*$. By (4.8), for $\varphi = s(0, q_0) \in c^*$, we have

$$\psi \in \widehat{D}^* P_{c_0}(x, h(x))(\varphi) \iff \limsup_{\substack{(u,v) \to (x, h(x)) \\ v \in P_{c_0}(u)}} \frac{\langle \psi, u-x \rangle - \langle s(0, q_0), v - h(x) \rangle}{\|u-x\| + \|v - h(x)\|} \leq 0$$

$$\iff \limsup_{\substack{(u,v) \to (x, h(x)) \\ v \in P_{c_0}(u)}} \frac{\langle \psi, u-x \rangle}{\|u-x\| + \|v - h(x)\|} \leq 0. \tag{4.10}$$

Case 1. There is a positive integer $m$ satisfying $p_m \neq 0$. Then, in the limit (4.10), we take a line segment direction as $u = x + ts(m, p_m)$, for $t \downarrow 0$. One has

$$u - x = ts(m, p_m) = h(x + ts(m, p_m)) - h(x) = ts(m, p_m) \to \theta, \text{ as } t \downarrow 0.$$

In the limit (4.10), we have

$$\limsup_{\substack{(u,v) \to (x, h(x)) \\ v \in P_{c_0}(u)}} \frac{\langle \psi, u-x \rangle}{\|u-x\| + \|v - h(x)\|}$$

$$\geq \limsup_{(u, h(u)) \to (x, h(x))} \frac{\langle \psi, u-x \rangle}{\|u-x\| + \|h(u) - h(x)\|}$$

$$= \limsup_{t \downarrow 0} \frac{\langle \psi, ts(m, p_m) \rangle}{\|ts(m, p_m)\| + \|ts(m, p_m)\|}$$

$$= \limsup_{t \downarrow 0} \frac{tp_m^2}{2t|p_m|}$$

$$= \frac{|p_m|}{2} > 0.$$

This contradicts to (4.10) for $\psi \in \widehat{D}^* P_{c_0}(x, h(x))(\varphi)$.

Case 2. Suppose $p_m = 0$, for all $m = 1, 2, \dots$. By the assumption that $\psi \neq \theta^*$, we have $p_0 \neq 0$. Then, in limit (4.10), we take a line segment direction as $u_t = x + \beta_{tp_0}$, for $t \downarrow 0$. One has

$$u_t - x = \beta_{tp_0} \to \theta, \text{ as } t \downarrow 0.$$

By Lemma 4.2, we have

$$h(u_t) = h(x + \beta_{tp_0}) = h(x) \quad \text{and} \quad h(u_t) - h(x) = \theta, \text{ for all } t > 0.$$

We calculate

$$\limsup_{\substack{(u,v)\to(x,h(x))\\v\in P_{c_0}(u)}} \frac{\langle\psi,u-x\rangle - \langle\varphi, v-h(x)\rangle}{\|u-x\| + \|v-h(x)\|}$$

$$\geq \limsup_{\substack{(u_t,v)\to(x,h(x))\\v\in P_{c_0}(u_t)}} \frac{\langle\psi,u_t-x\rangle - \langle\varphi, v-h(x)\rangle}{\|u_t-x\| + \|v-h(x)\|}$$

$$= \limsup_{(x+\beta_{tp_0},h(x+\beta_{tp_0}))\to(x,h(x))} \frac{\langle\psi,x+\beta_{tp_0}-x\rangle - \langle\varphi, h(x+\beta_{tp_0})-h(x)\rangle}{\|x+\beta_{tp_0}-x\| + \|h(x+\beta_{tp_0})-h(x)\|}$$

$$= \limsup_{t\downarrow 0} \frac{\langle\psi, \beta_{tp_0}\rangle - \langle\varphi, h(x)-h(x)\rangle}{\|\beta_{tp_0}\| + \|h(x)-h(x)\|}$$

$$= \limsup_{t\downarrow 0} \frac{\langle\psi, \beta_{tp_0}\rangle}{\|\beta_{tp_0}\|}$$

$$= \limsup_{t\downarrow 0} \frac{p_0 \lim_{n\to\infty} tp_0 + \sum_{n=1}^{\infty} p_n tp_0}{t|p_0|}$$

$$= \frac{p_0 \lim_{n\to\infty} p_0}{|p_0|}$$

$$= |p_0| > 0.$$

This contradicts to (4.10). By case 1 and case 2, we obtain

$$\psi \neq \theta^* \implies \psi \notin \widehat{D}^*P_{c_0}(x,h(x))(\varphi), \text{ for } \varphi = s(0,q_0) \in c^* \text{ with } q_0 \in \mathbb{R}. \quad (4.11)$$

By (4.9) and (4.11), (i) is proved:

$$\widehat{D}^*P_{c_0}(x,h(x))(\varphi) = \{\theta^*\}, \text{ for } \varphi = s(0,q_0) \in c^* \text{ with } q_0 \in \mathbb{R}.$$

Proof of (ii). Let $x \in c$. For any $\varphi = (q_0, q_1, q_2, \ldots) \in c^*$ with $\varphi \neq s(0,q_0)$, for any $q_0 \in \mathbb{R}$, which includes that $\varphi \neq \theta^* = s(0,0)$. Then, there is at least one positive integer $m$ with $q_m \neq 0$. At first, for such $\varphi$, we show

$$\theta^* \notin \widehat{D}^*P_{c_0}(x,h(x))(\varphi). \quad (4.12)$$

By definition, we have

$$\theta^* \in \widehat{D}^*P_{c_0}(x,h(x))(\varphi) \iff \limsup_{\substack{(u,v)\to(x,h(x))\\v\in P_{c_0}(u)}} \frac{-\langle\varphi, v-h(x)\rangle}{\|u-x\| + \|v-h(x)\|} \leq 0. \quad (4.13)$$

Then, in limit (4.13), we take a line segment direction as $u_t = x + ts(m, -q_m)$, for $t \downarrow 0$. One has

$$u_t - x = h(x + ts(m, -q_m)) - h(x) = ts(m, -q_m) \to \theta, \text{ as } t \downarrow 0.$$

We have

$$\limsup_{\substack{(u,v)\to(x,h(x))\\v\in P_{c_0}(u)}} \frac{-\langle \varphi, v-h(x)\rangle}{\|u-x\|+\|v-h(x)\|}$$

$$\geq \limsup_{(u_t,h(u_t))\to(x,h(x))} \frac{-\langle \varphi, h(u_t)-h(x)\rangle}{\|u_t-x\|+\|h(u_t)-h(x)\|}$$

$$\geq \limsup_{t\downarrow 0} \frac{-\langle \varphi, ts(m,-q_m)\rangle}{\|ts(m,-q_m)\|+\|ts(m,-q_m)\|}$$

$$= \limsup_{t\downarrow 0} \frac{tq_m^2}{2t|q_m|}$$

$$= \frac{|q_m|}{2} > 0.$$

This contradicts to (4.13), which proves (4.12).

Next, under condition that $\varphi = (q_0, q_1, q_2, ...) \in c^*$ with $\varphi \neq s(0, q_0)$, for any $q_0 \in \mathbb{R}$, we show

$$\psi \notin \widehat{D}^* P_{c_0}(x, h(x))(\varphi), \text{ for any } \psi \in c^*\setminus\{\theta^*\}. \tag{4.12.1}$$

Let $\psi = (p_0, p_1, p_2, ...) \in c^*\setminus\{\theta^*\}$. The proof is divided into the following two cases.

Case 1. $\sum_{n=0}^{\infty} p_n \neq 0$. In this case, we prove the following result stranger than (4.12.1).

$$\sum_{n=0}^{\infty} p_n \neq 0 \implies \psi \notin \widehat{D}^* P_{c_0}(x, h(x))(\varphi), \text{ for any } x \in c \text{ and for any } \varphi \in c^*. \tag{4.14}$$

Let $p := \sum_{n=0}^{\infty} p_n$. Let $x \in c$ and $\varphi \in c^*$ be arbitrarily given. By definition, we have

$$\psi \in \widehat{D}^* P_{c_0}(x, h(x))(\varphi) \iff \limsup_{\substack{(u,v)\to(x,h(x))\\v\in P_{c_0}(u)}} \frac{\langle \psi, u-x\rangle - \langle \varphi, v-h(x)\rangle}{\|u-x\|+\|v-h(x)\|} \leq 0. \tag{4.15}$$

Then, in limit (4.15), we take a line segment direction as $u_t = x + \beta_{tp}$, for $t \downarrow 0$. One has

$$u_t - x = \beta_{tp} \to \theta, \text{ as } t \downarrow 0.$$

By Lemma 4.2, we have

$$h(u_t) = h(x + \beta_{tp}) = h(x) \quad \text{and} \quad h(u_t) - h(x) = \theta, \text{ for all } t > 0.$$

We calculate

$$\limsup_{\substack{(u,v)\to(x,h(x))\\v\in P_{c_0}(u)}} \frac{\langle \psi, u-x\rangle - \langle \varphi, v-h(x)\rangle}{\|u-x\|+\|v-h(x)\|}$$

$$\geq \limsup_{\substack{(u_t,v)\to(x,h(x))\\v\in P_{c_0}(u_t)}} \frac{\langle \psi, u_t-x\rangle - \langle \varphi, v-h(x)\rangle}{\|u_t-x\|+\|v-h(x)\|}$$

$$= \limsup_{(x+\beta_{tp},h(x+\beta_{tp}))\to(x,h(x))} \frac{\langle \psi, x+\beta_{tp}-x\rangle - \langle \varphi, h(x+\beta_{tp})-h(x)\rangle}{\|x+\beta_{tp}-x\| + \|h(x+\beta_{tp})-h(x)\|}$$

$$= \limsup_{t\downarrow 0} \frac{\langle \psi, \beta_{tp}\rangle - \langle \varphi, h(x)-h(x)\rangle}{\|\beta_{tp}\| + \|h(x)-h(x)\|}$$

$$= \limsup_{t\downarrow 0} \frac{\langle \psi, \beta_{tp}\rangle}{\|\beta_{tp}\|}$$

$$= \limsup_{t\downarrow 0} \frac{p_0 \lim_{n\to\infty} tp + \sum_{n=1}^{\infty} p_n tp}{t|p|}$$

$$= \frac{p \sum_{n=0}^{\infty} p_n}{|p|}$$

$$= |p| > 0.$$

By (4.15), this implies $\psi \notin \widehat{D}^* P_{C_0}(x, h(x))(\varphi)$, which proves (4.14).

Case 2. $\sum_{n=0}^{\infty} p_n = 0$. In this case, by (4.12), for such a $\psi = (p_0, p_1, p_2, ...) \ c^*\backslash\{\theta^*\}$, we prove

$$\sum_{n=0}^{\infty} p_n = 0 \implies \psi \notin \widehat{D}^* P_{C_0}(x, h(x))(\varphi), \text{ for } \varphi \neq s(0, q_0), \text{ for any } q_0 \in \mathbb{R}. \quad (4.16)$$

For any positive integer $m$, in limit (4.15), we take a line segment direction as $u_t = x + ts(m, 1)$, for $t \downarrow 0$. One has

$$u_t - x = h(x + ts(m, 1)) - h(x) = ts(m, 1) \to \theta, \text{ as } t \downarrow 0.$$

Then,

$$\psi \in \widehat{D}^* P_{C_0}(x, h(x))(\varphi) \implies 0 \geq \limsup_{\substack{(u,v)\to(x,h(x))\\ v\in P_{C_0}(u)}} \frac{\langle \psi, u-x\rangle - \langle \varphi, v-h(x)\rangle}{\|u-x\| + \|v-h(x)\|}$$

$$\geq \limsup_{\substack{(u_t,h(u_t))\to(x,h(x))\\ v\in P_{C_0}(u)}} \frac{\langle \psi, u_t-x\rangle - \langle \varphi, h(u_t)-h(x)\rangle}{\|u_t-x\| + \|h(u_t)-h(x)\|}$$

$$= \limsup_{t\downarrow 0} \frac{\langle \psi, ts(m,1)\rangle - \langle \varphi, ts(m,1)\rangle}{\|ts(m,1)\| + \|ts(m,1)\|}$$

$$= \limsup_{t\downarrow 0} \frac{t(p_m - q_m)}{2t}$$

$$= \frac{p_m - q_m}{2}.$$

This implies that

$$\psi \in \widehat{D}^* P_{C_0}(x, h(x))(\varphi) \implies p_m \leq q_m, \text{ for } m = 1, 2, ... \quad (4.17)$$

For any positive integer $m$, in limit (4.15), we take a line segment direction as $u_t = x + ts(m, -1)$, for $t \downarrow 0$. One has

$$u_t - x = h(x + ts(m, -1)) - h(x) = ts(m, -1) \to \theta, \text{ as } t \downarrow 0.$$

Then,

$$\psi \in \widehat{D}^* P_{c_0}(x, h(x))(\varphi) \implies 0 \geq \limsup_{\substack{(u,v) \to (x, h(x)) \\ v \in P_{c_0}(u)}} \frac{\langle \psi, u-x \rangle - \langle \varphi, v - h(x) \rangle}{\|u-x\| + \|v - h(x)\|}$$

$$\geq \limsup_{t \downarrow 0} \frac{\langle \psi, ts(m,-1) \rangle - \langle \varphi, ts(m,-1) \rangle}{\| ts(m,1) \| + \| ts(m,1) \|}$$

$$= \limsup_{t \downarrow 0} \frac{t(-p_m + q_m)}{2t}$$

$$= \frac{-p_m + q_m}{2}.$$

This implies that

$$\psi \in \widehat{D}^* P_{c_0}(x, h(x))(\varphi) \implies p_m \geq q_m, \text{ for } m = 1, 2, \ldots \tag{4.18}$$

By (4.17) and (4.18), we have

$$\psi \in \widehat{D}^* P_{c_0}(x, h(x))(\varphi) \implies p_m = q_m, \text{ for } m = 1, 2, \ldots \tag{4.19}$$

By (4.19), for any $y \in c_0$, one has,

$$\psi \in \widehat{D}^* P_{c_0}(x, h(x))(\varphi) \implies \langle \varphi, y \rangle - \langle \psi, y \rangle = q_0 L(y) - p_0 L(y) = 0.$$

In particular, for any $x, u \in c$, we have

$$\langle \varphi, v - h(x) \rangle = \langle \psi, v - h(x) \rangle, \text{ for any } v \in P_{c_0}(u). \tag{4.20}$$

By $\sum_{n=0}^{\infty} p_n = 0$ and by (4.19) and (4.20), this implies

$$\psi \in \widehat{D}^* P_{c_0}(x, h(x))(\varphi) \implies 0 \geq \limsup_{\substack{(u,v) \to (x, h(x)) \\ v \in P_{c_0}(u)}} \frac{\langle \psi, u-x \rangle - \langle \varphi, v - h(x) \rangle}{\|u-x\| + \|v - h(x)\|}$$

$$= \limsup_{\substack{(u,v) \to (x, h(x)) \\ v \in P_{c_0}(u)}} \frac{\langle \psi, u-x \rangle - \langle \psi, v - h(x) \rangle}{\|u-x\| + \|v - h(x)\|}$$

$$= \limsup_{\substack{(u,v) \to (x, h(x)) \\ v \in P_{c_0}(u)}} \frac{\langle \psi, u-v-(x-h(x)) \rangle}{\|u-x\| + \|v - h(x)\|}$$

$$
\begin{aligned}
&= \limsup_{\substack{(u,v)\to(x,h(x))\\ v\in P_{C_0}(u)}} \frac{\langle \psi, u-v-\beta_{L(x)}\rangle}{\|u-x\|+\|v-h(x)\|}\\
&= \limsup_{\substack{(u,v)\to(x,h(x))\\ v\in P_{C_0}(u)}} \frac{\langle \psi, u-v\rangle - L(x)\sum_{n=0}^{\infty} p_n}{\|u-x\|+\|v-h(x)\|}\\
&= \limsup_{\substack{(u,v)\to(x,h(x))\\ v\in P_{C_0}(u)}} \frac{\langle \psi, u-v\rangle}{\|u-x\|+\|v-h(x)\|}.
\end{aligned}
$$

This is,

$$\psi \in \widehat{D}^*P_{C_0}(x,h(x))(\varphi) \implies \limsup_{\substack{(u,v)\to(x,h(x))\\ v\in P_{C_0}(u)}} \frac{\langle \psi, u-x\rangle}{\|u-x\|+\|v-h(x)\|} \leq 0. \quad (4.21)$$

By (4.21), similar to the proof of (4.18), one has

$$\psi \in \widehat{D}^*P_{C_0}(x,h(x))(\varphi) \implies p_m = 0, \text{ for } m = 1, 2, \ldots \quad (4.22)$$

By the assumption that $\sum_{n=0}^{\infty} p_n = 0$, (4.22) implies that $p_0 = 0$, and therefore, $\psi = \theta^*$. This contradicts to the assumption that $\psi = (p_0, p_1, p_2, \ldots) \, c^*\setminus\{\theta^*\}$, which proves (4.16). By (4.12), (4.14) and (4.16), we proved that, for any $x \in c$,

$$\widehat{D}^*P_{C_0}(x,h(x))(\varphi) = \emptyset, \text{ for any } \varphi \in c^* \text{ with } \varphi \neq s(0,p_0), \text{ for } p_0 \in \mathbb{R}. \qquad \square$$

## 5. Banach space $C[0, 1]$

### 5.1. On the Chebyshev's Equioscillation Theorem

Let $(C[0, 1], \|\cdot\|)$ be the Banach space of all continuous real valued functions on $[0, 1]$ with respect to the standard Boral $\sigma$-field $\Sigma$ and with the maximum norm

$$\|f\| = \max_{0\leq t \leq 1} |f(t)|, \text{ for any } f \in C[0, 1].$$

The dual space of $C[0, 1]$ is denoted by $C^*[0, 1]$ that is $rca[0, 1]$, in which the considered $\sigma$-field is the standard Boral $\sigma$-field $\Sigma$ on $[0, 1]$ including all closed and open subsets of $[0, 1]$. Let $\langle \cdot, \cdot \rangle$ denote the real canonical pairing between $C^*[0, 1]$ and $C[0, 1]$. By the Riesz Representation Theorem (see Theorem IV.6.3 in Dunford and Schwartz [8]), for any $\varphi \in C^*[0, 1]$, there is a real valued, regular and countable additive functional $\mu \in rca[0, 1]$, which is defined on the given $\sigma$-field $\Sigma$ on $[0, 1]$, such that

$$\langle \varphi, f \rangle = \int_0^1 f(t)\mu(dt), \text{ for any } f \in C[0, 1]. \quad (5.1)$$

The norm on $C^*[0, 1]$ is denoted by $\|\cdot\|_*$ satisfying that, for any $\varphi \in C^*$, the norm of $\varphi$ in $C^*$ is

$$\|\varphi\|_* = \|\mu\|_* := v(\mu, [0, 1]). \quad (5.2)$$

Where, $v(\mu, [0, 1])$ is the total variation of $\mu$ on $[0, 1]$, which is defined by (see page 160 in [8])

$$v(\mu, [0, 1]) = \sup_{E \in \Sigma} |\mu(E)|.$$

Throughout this subsection, without any special mention, we shall identify $\varphi$ and $\mu$ in (5.1). We say that $\mu \in C^*[0, 1]$ ((it is $rca[0, 1]$) and (5.1) is rewritten as

$$\langle \mu, f \rangle = \int_0^1 f(t)\mu(dt), \text{ for any } f \in C[0, 1]. \tag{5.1}$$

The origin of $C[0, 1]$ is denoted by $\theta$, which is the constant function defined on $[0, 1]$ with value 0. The origin of the dual space $C^*$ is denoted by $\theta^*$, which is also a constant functional on $\Sigma$ with value 0. This is,

$$\theta^*(E) = 0, \text{ for every } E \in \Sigma.$$

For any $f \in C[0, 1]$, define

$$M(f) = \{t \in [0, 1]: |f(t)| = \|f\|\}.$$

By the continuity of $f$ on $[0, 1]$, $M(f)$ is a nonempty closed subset of $[0, 1]$, which is called the maximizing set of $f$. The following lemma shows the connections between the normalized duality mapping $J: C \to 2^{C^*} \setminus \{\emptyset\}$ and the maximizing sets.

**Lemma 5.1**. *For any $f \in C[0, 1]$ with $\|f\| > 0$, then*

(a) $\mu \in J(f) \implies v(\mu, [0, 1] \setminus M(f)) = 0$;
(b) *For any $\{t_j \in M(f): j = 1, 2, \ldots, m\} \subseteq M(f)$, for some positive integer $m$, define $\mu \in rca[0, 1]$ by*

$$\mu(t_j) = \alpha_j f(t_j), \text{ for } j = 1, 2, \ldots, m \quad \text{and} \quad \mu([0, 1] \setminus \{t_j: j = 1, 2, \ldots, m\}) = 0,$$

*where $\alpha_j > 0$, for $j = 1, 2, \ldots, m$ satisfying $\sum_{j=1}^m \alpha_j = 1$, then, $\mu \in J(f)$.*

(c) *If $M(f)$ is not a singleton, then $J(f)$ is an infinite set.*

*Proof*. Proof of (a). Assume, by the way of contradiction, that $\mu \in J(f)$ and $v(\mu, [0, 1] \setminus M(f)) > 0$. Then, by (5.1) and (5.2), we have

$$\|f\|^2 = \langle \mu, f \rangle = \|\mu\|_*^2 = (v(\mu, [0, 1]))^2.$$

This implies that

$$\|f\|^2 = \langle \mu, f \rangle$$
$$= \int_0^1 f(t)\mu(dt)$$
$$= \int_{M(f)} f(t)\mu(dt) + \int_{[0,1]\setminus M(f)} f(t)\mu(dt)$$
$$\leq \|f\| v(\mu, M(f)) + \int_{[0,1]\setminus M(f)} f(t)\mu(dt)$$

$$< \|f\|v(\mu, M(f)) + \|f\|v(\mu, [0, 1]\backslash M(f))$$
$$= \|f\|v(\mu, [0, 1])$$
$$= \|f\|^2.$$

It is a contradiction.

Proof of (b). By definitions of $\mu$ and $\|\mu\|_*$, we calculate

$$\|\mu\|_* = v(\mu, [0, 1]) = \sum_{j=1}^m |\mu(t_j)| = \sum_{j=1}^m |\alpha_j f(t_j)| = \sum_{j=1}^m \alpha_j \|f\| = \|f\|. \tag{5.3}$$

On the other hand, we have

$$\langle \mu, f \rangle = \sum_{j=1}^m f(t_j)\mu(t_j) = \sum_{j=1}^m f(t_j)\alpha_j f(t_j) = \sum_{j=1}^m \alpha_j \|f\|^2 = \|f\|^2. \tag{5.4}$$

By (5.3) and (5.4), it proves $\mu \in J(f)$. Then, part (c) follows from part (b) immediately. □

For any given nonnegative integer $n$, let $\mathcal{P}_n$ denote the closed subspace of $C[0, 1]$ that consists of all real coefficients' polynomials of degree less than or equal to $n$. Let $P_{\mathcal{P}_n}$ denote the metric projection from $C[0, 1]$ to $\mathcal{P}_n$. Let $f \in C[0, 1]$ and $p \in \mathcal{P}_n$. Then $p \in P_{\mathcal{P}_n}(f)$ if and only if

$$\|f - p\| = \min\{\|f - q\|: q \in \mathcal{P}_n\}.$$

That is, if $p \in P_{\mathcal{P}_n}(f)$, then $p$ is the polynomial of best uniform approximation for $f$ satisfying,

$$\|f - p\|$$
$$= \max_{0 \le t \le 1} |f(t) - p(t)|$$
$$= \min\{\|f - q\|: q(t) = \sum_{k=0}^n c_k t^k \in \mathcal{P}_n\}$$
$$= \min\{\max_{0 \le t \le 1} |f(t) - \sum_{k=0}^n c_k t^k|: \sum_{k=0}^n c_k t^k \in \mathcal{P}_n\}.$$

The metric projection operator $P_{\mathcal{P}_n}$ from $C[0, 1]$ to $\mathcal{P}_n$ has the following existence properties.

**Proposition 5.2**. *Let n be a nonnegative integer. Then*

(a) $P_{\mathcal{P}_n}(p) = p$, *for any* $p \in \mathcal{P}_n$;
(b) $P_{\mathcal{P}_n}(f) \ne \emptyset$, *for any* $f \in C[0, 1]$.

*Proof.* Part (a) is clear. We only prove part (b) for $f \in C[0, 1]\backslash\mathcal{P}_n$. Since $\mathcal{P}_n$ is a closed subspace of $C[0, 1]$, then, for any $f \in C[0, 1]\backslash\mathcal{P}_n$, we have

$$A := \min\{\|f - q\|: q \in \mathcal{P}_n\} > 0.$$

We take a sequence $\{q_m\} \subseteq \mathcal{P}_n$ such that

$$\|f - q_m\| \downarrow A, \text{ as } m \to \infty.$$

This implies that $\{\|q_m\|\}$ is bounded. Since every $q'_m$ is a polynormal with degree less than or equal to $n - 1$, it follows that $\{\|q'_m\|\}$ is bounded. Let $d = \max\{\|q'_m\|: m = 1, 2, 3 ...\}$. For any given $\varepsilon > 0$, we take $\delta > 0$ with $\delta < \frac{\varepsilon}{d}$. By the mean value theorem, this implies that, for any $s, t \in [0, 1]$, we have

$$|s - t| < \delta \implies |q_m(s) - q_m(t)| < \varepsilon, \text{ for } m = 1, 2, \ldots.$$

This implies that $\{\|q_m\|\}$ is a bounded and equicontinuous subset in $C[0, 1]$. By Azsela-Ascali Theorem, $\{\|q_m\|\}$ is conditionally compact. Hence, there is a subsequence $\{q_{m_i}\} \subseteq \{q_m\} \subseteq \mathcal{P}_n$ and there is $q \in \mathcal{P}_n$ such that $q_{m_i} \to q$, as $i \to \infty$. We have

$$\|f - q\| \leq \|f - q_{m_i}\| + \|q - q_{m_i}\| \to A, \text{ as } i \to \infty.$$

This implies $\|f - q\| = A$. That is $q \in P_{\mathcal{P}_n}(f)$. □

In 1859, Chebyshev proved the well-known theorem called the Chebyshev's Equioscillation theorem. This theorem provides the criterion for the necessary and sufficient conditions for the solutions of best approximation of continuous functions by polynomials with degree less than or equal to $n$. We list the revised version of the Chebyshev's Equioscillation Theorem with respect to the metric projection $P_{\mathcal{P}_n}$ (see [5], [6], [29]).

**Chebyshev's Equioscillation Theorem**. *Let $f(t)$ be a continuous function on $[0,1]$ and, for $p(t) = \sum_{k=0}^{n} a_k t^k \in \mathcal{P}_n$, suppose that*

$$A(f, p) := \|f - p\| = \max_{0 \leq t \leq 1} |f(t) - p(t)|.$$

*Then, $p \in P_{\mathcal{P}_n}(f)$ if and only if, there are $n+2$ points $0 \leq t_0 < t_1 < \cdots < t_{n+1} \leq 1$ such that*

$$f(t_i) - p(t_i) = \epsilon A(f, p)(-1)^i, i = 0, 1, 2, \ldots, n+1, \tag{5.5}$$

*where $\epsilon = 1$ or $-1$.*

The set $\{t_0, t_1, \ldots, t_{n+1}\}$ satisfying (5.5) is called a $n$-Chebyshev set of $f$ with respect to $p$, which is denoted by $S(f, p)$. It is clear that

$$S(f, p) \subseteq M(f - p), \text{ for } p \in P_{\mathcal{P}_n}(f).$$

In the following example, we show that $S(f, p)$ may not be unique, in general.

**Example 5.3**. Let $n = 1$. Let $f(t) = \sin(4\pi t) \in C[0, 1]$. Then, $P_{\mathcal{P}_1}(f) = \theta$ with $A(f, \theta) = 1$. One has that there are two 1-Chebyshev sets of $f$ with respect to $p$, which are

$$\left\{\frac{\pi}{2}, \frac{3\pi}{2}, \frac{5\pi}{2}\right\} \text{ with } \epsilon = 1 \quad \text{and} \quad \left\{\frac{3\pi}{2}, \frac{5\pi}{2}, \frac{7\pi}{2}\right\} \text{ with } \epsilon = -1.$$

Since the Banach space $C[0, 1]$ is neither reflexive, nor strictly convex and smooth, so, it is not for sure that the metric projection $P_{\mathcal{P}_n}$ from $C[0, 1]$ to $\mathcal{P}_n$ is a single-valued mapping or not. However, by Proposition 5.2, the following theorem shows that $\mathcal{P}_n$ is indeed a single-valued mapping (See Bray B. R., Minimax approximation theory, Nov. 27, 2016). We provide a proof of the following theorem in the appendix of this paper.

**Theorem 5.4**. *Let $n$ be a nonnegative integer. The metric projection $P_{\mathcal{P}_n}: C[0, 1] \to \mathcal{P}_n$ is a single-valued mapping.*

By the Chebyshev's Equioscillation Theorem, we have the following approximation properties.

**Proposition 5.5.** *Let n be a nonnegative integer. Let $f \in C[0, 1]$ with $p = P_{\mathcal{P}_n}(f)$. Then,*

(i)
$$\beta p = P_{\mathcal{P}_n}(\beta f), \text{ for any } \beta \in \mathbb{R},$$

$$A(\beta f, \beta p) = |\beta| A(f, p) \quad \text{and} \quad S(\beta f, \beta p) = S(f, p),$$

$$\epsilon(\beta f, \beta p) = \epsilon(f, p), \text{ if } \beta \geq 0 \quad \text{and} \quad \epsilon(\beta f, \beta p) = -\epsilon(f, p), \text{ if } \beta < 0.$$

(ii)
$$p + q = P_{\mathcal{P}_n}(f + q), \text{ for any } q \in \mathcal{P}_n,$$

$$A(f + q, p + q) = A(f, p) \quad \text{and} \quad S(f + q, p + q) = S(f, p).$$

*In particular, we have*

$$p + \lambda = P_{\mathcal{P}_n}(f + \lambda), \quad \text{for any } \lambda \in \mathbb{R};$$

(iii) *Let $f \in C[0, 1] \setminus \mathcal{P}_n$ and $p \in \mathcal{P}_n$ with $p = P_{\mathcal{P}_n}(f)$. Then*

$$p = P_{\mathcal{P}_n}(\alpha f + (1 - \alpha)p), \text{ for any } \alpha \in (0, 1).$$

*Proof.* The proofs of parts (i) and (ii) are straightforward and they are omitted here. We only show part (iii) by using parts (i) and (ii).

For any given $\alpha \in (0, 1)$, by (i), we have $\alpha p = P_{\mathcal{P}_n}(\alpha f)$. Then, by (ii), we have

$$p = \alpha p + (1 - \alpha)p = P_{\mathcal{P}_n}(\alpha f + (1 - \alpha)p). \qquad \square$$

Let $n$ be an arbitrarily given nonnegative integer. For $f \in C[0, 1]$ and $p \in \mathcal{P}_n$ with $p = P_{\mathcal{P}_n}(f)$. Let $S(f, p) = \{t_i \in [0,1]: i = 0, 1, 2, \ldots, n + 1\} \subseteq M(f - p)$. We define a regular countable additive functional $\mu_p^f \in rca[0, 1]$ by

$$\mu_p^f(t_i) = \frac{1}{n+2}(f(t_i) - p(t_i)), \text{ for } i = 0,1, 2, \ldots, n+1 \text{ and } \mu_p^f([0, 1] \setminus S(f,p)) = 0.$$

By part (b) of Lemma 5.1, we have

$$\mu_p^f \in J(f-p), \text{ for } p = P_{\mathcal{P}_n}(f).$$

By Chebyshev's Equioscillation Theorem, we have

$$\mu_p^f(t_i) = \frac{\epsilon}{n+2} A(f, p)(-1)^i = \frac{\epsilon}{n+2} \|f - p\|(-1)^i, \text{ for } i = 0, 1, 2, \ldots, n+1.$$

We provide some examples below to practice the Chebyshev's Equioscillation Theorem.

**Example 5.6.** Let $n = 2$. Let $f(t) = t^3 - t \in C[0, 1] \setminus \mathcal{P}_2$. One sees that $a = \frac{2}{3}$ is a real solution of the following equation:

$$\frac{a^3}{54} - \frac{(2-a)^2}{4} = 0. \tag{5.6}$$

Then, $p := P_{\mathcal{P}_2}(f)$ is given by

$$p(t) = at^2 - \left(1 + \frac{a^2}{4}\right)t + \frac{(2-a)^2}{8} = \frac{2}{3}t^2 - \frac{25}{16}t + \frac{1}{32}, \text{ for all } t \in [0, 1].$$

$A(f, p) = \|f - p\| = \frac{(2-a)^2}{8} = \frac{1}{32}$. There is a unique 2-Chebyshev set of $f$ with respect to $p$, which is $S(f, p) = \left\{0, \frac{a}{6}, \frac{a}{2}, 1\right\} = \left\{0, \frac{1}{4}, \frac{3}{4}, 1\right\}$ with $\epsilon = -1$.

*Proof.* Notice that the polynormal $\frac{a^3}{54} - \frac{(2-a)^2}{4}$ in (5.6) has negative value at point $t = 1$ and has positive value at point $t = 2$. This implies that the equation (5.6) has a solution in $(1, 2)$, which ensures the existence of the desired $a = \frac{2}{3} \in (1, 2)$ satisfying (5.6).

Define $g := f - p \in C[0, 1]$ by

$$g(t) = t^3 - t - \left(at^2 - \left(1 + \frac{a^2}{4}\right)t + \frac{(2-a)^2}{8}\right)$$

$$= t^3 - at^2 + \frac{a^2}{4}t - \frac{(2-a)^2}{8}, \text{ for all } t \in [0, 1].$$

$g$ has two critical points $t_1 = \frac{a}{6}$ and $t_2 = \frac{a}{2}$. By $a \in (1, 2)$, both $t_1$ and $t_2$ are in $(0,1)$. By (5.6), the extreme value of $g$ at the first critical point $t_1$ is

$$g\left(\frac{a}{6}\right) = \left(\frac{a}{6}\right)^3 - a\left(\frac{a}{6}\right)^2 + \frac{a^2}{4}\frac{a}{6} - \frac{(2-a)^2}{8}$$

$$= \frac{a^3}{54} - \frac{(2-a)^2}{8}$$

$$= \frac{a^3}{54} - \frac{(2-a)^2}{4} + \frac{(2-a)^2}{8}$$

$$= \frac{(2-a)^2}{8}. \tag{5.7}$$

The extreme value of $g$ at the second critical point $t_2$ is

$$g\left(\frac{a}{2}\right) = \left(\frac{a}{2}\right)^3 - a\left(\frac{a}{2}\right)^2 + \frac{a^2}{4}\frac{a}{2} - \frac{(2-a)^2}{8} = -\frac{(2-a)^2}{8}. \tag{5.8}$$

The value of $g$ at the left end point of $[0, 1]$ is

$$g(0) = -\frac{(2-a)^2}{8}. \tag{5.9}$$

The value of $g$ at the right end point of $[0, 1]$ is

$$g(1) = 1 - a + \frac{a^2}{4} - \frac{(2-a)^2}{8} = \frac{(2-a)^2}{8}. \tag{5.10}$$

By (5.7)–(5.10), we obtain

$$A(f, p) = \|f - p\| = \|g\| = \frac{(2-a)^2}{8} = \frac{1}{32}.$$

By, (5.9), (5.7), (5.8) and (5.10), we have a 2-Chebyshev set of $f$ with respect to $p$, which is $S(f, p) = \left\{0, \frac{a}{6}, \frac{a}{2}, 1\right\} = \left\{0, \frac{1}{4}, \frac{3}{4}, 1\right\}$ with $\epsilon = -1$. □

As an application of Proposition 5.2, we create the following example by using the previous example.

**Example 5.7.** Let $n = 2$. Let $f$ and $p := P_{\mathcal{P}_2}(f)$ be given in Example 5.6. Let

$$h(t) = f(t) + t = t^3 \in C[0, 1]\backslash \mathcal{P}_2.$$

Let $q := P_{\mathcal{P}_2}(h)$. By Proposition 5.2 and Example 5.6, we have

$$q(t) = p(t) + t = at^2 - \frac{a^2}{4}t + \frac{(2-a)^2}{8} = \frac{2}{3}t^2 - \frac{9}{16}t + \frac{1}{32}, \text{ for all } t \in [0, 1].$$

We also have the 2-Chebyshev set of $h$ with respect to $q$ is $S(h, q) = S(f, p) = \left\{0, \frac{a}{6}, \frac{a}{2}, 1\right\} = \left\{0, \frac{1}{4}, \frac{3}{4}, 1\right\}$ with $\epsilon = -1$ and

$$A(h, q) = \|h - q\| = A(f, p) = \|f - p\| = \frac{(2-a)^2}{8} = \frac{1}{32}.$$

**Lemma 5.8.** *Let $n$ be any positive integer. For any $q \in \mathcal{P}_n$, there is $g \in C[0, 1]\backslash \mathcal{P}_n$ such that $q = P_{\mathcal{P}_n}(g)$.*

*Proof.* For any given positive integer $n$ and for any $q \in \mathcal{P}_n$, we take a positive integer $m > n$. Let $f(t) = \sin(2m\pi t)$. We have $\theta = P_{\mathcal{P}_n}(f)$. Let $g = f + q$. We have $g \in C[0, 1]\backslash \mathcal{P}_n$. By Proposition 5.5, we have

$$q = \theta + q = P_{\mathcal{P}_n}(f + q) = P_{\mathcal{P}_n}(g). \qquad \square$$

Consequently, from Lemma 5.8, we obtain the following results.

**Corollary 5.9.** *Let $n$ be a positive integer. The interior of $\mathcal{P}_n$ is empty. That is,*

$$(\mathcal{P}_n)^o = \emptyset, \text{ for any positive integer } n.$$

*Proof.* Assume that there is $q \in \mathcal{P}_n$ and $\varepsilon > 0$ such that

$$\{h \in C[0, 1]: \|q - h\| < \varepsilon\} \subseteq \mathcal{P}_n. \tag{5.11}$$

By Lemma 5.8, there is $g \in C[0, 1]\backslash \mathcal{P}_n$ such that $q = P_{\mathcal{P}_n}(g)$. By part (iii) of Proposition 5.5,

$$q = P_{\mathcal{P}_n}(\alpha g + (1 - \alpha)q), \text{ for any } \alpha \in (0, 1).$$

We take $\alpha$ with $0 < \alpha < \frac{\varepsilon}{\|q-g\|}$. Then $\|q - (\alpha g + (1 - \alpha)q)\| = \alpha\|q - g\| < \varepsilon$. By (5.11),

$$\alpha g + (1 - \alpha)q \in \mathcal{P}_n, \text{ for any } \alpha \text{ with } 0 < \alpha < \frac{\varepsilon}{\|q-g\|}.$$

This contradicts to the assumption that $g \in C[0, 1]\backslash \mathcal{P}_n$. $\qquad \square$

### 5.2. The Mordukhovich and Gâteaux directional derivatives of the metric projection in $C[0, 1]$

In this subsection, we use the properties of the single-valued metric projection $P_{\mathcal{P}_n}: C[0, 1] \to \mathcal{P}_n$ studied in the previous subsection to investigate the solutions of its Mordukhovich derivatives and its Gâteaux directional derivatives.

**Theorem 5.10**. *Let n be a positive integer. Let $f \in C[0, 1]$ and $p \in \mathcal{P}_n$ with $p = P_{\mathcal{P}_n}(f)$. For any $\mu, \gamma \in C^*[0, 1]$, we have*

(i) $\mu([0,1]) \neq \gamma([0,1]) \implies \gamma \notin \widehat{D}^* P_{\mathcal{P}_n}(f,p)(\mu)$ *and* $\mu \notin \widehat{D}^* P_{\mathcal{P}_n}(f,p)(\gamma)$;

(ii) $\mu([0,1]) \neq 0 \implies \theta^* \notin \widehat{D}^* P_{\mathcal{P}_n}(f,p)(\mu)$ *and* $\mu \notin \widehat{D}^* P_{\mathcal{P}_n}(f,p)(\theta^*)$;

(iii) $\langle \gamma, f - p \rangle < 0 \implies \gamma \notin \widehat{D}^* P_{\mathcal{P}_n}(f,p)(\mu)$, *for any $\mu \in C^*$*.

*Proof.* Proof of (i). For any real number $\lambda$, we consider it as a constant function defined on $[0, 1]$ with value $\lambda$. Then, for any $\mu \in C^*$, we have

$$\langle \mu, \lambda \rangle = \lambda \mu([0,1]), \text{ for any } \lambda \in \mathbb{R}. \tag{5.12}$$

For any $\mu, \gamma \in C^*$, by definition, we have

$$\gamma \in \widehat{D}^* P_{\mathcal{P}_n}(f,p)(\mu) \iff \limsup_{\substack{(g,q) \to (f,p) \\ q = P_{\mathcal{P}_n}(g)}} \frac{\langle \gamma, g-f \rangle - \langle \mu, q-p \rangle}{\|g-f\| + \|q-p\|} \leq 0. \tag{5.13}$$

Suppose that $\mu([0,1]) < \gamma([0,1])$. Then, in limit (5.13), we take a line segment direction as $g_\lambda = f + \lambda$, for $\lambda \downarrow 0$. One has

$$g_\lambda - f = \lambda \to \theta, \text{ as } \lambda \downarrow 0.$$

By Proposition 5.5, we have

$$p + \lambda = P_{\mathcal{P}_n}(f + \lambda) = P_{\mathcal{P}_n}(g_\lambda), \text{ for any } \lambda > 0.$$

It is clear that $p + \lambda \to p$, as $\lambda \downarrow 0$. By (5.12), we calculate the limit in (5.13),

$$\limsup_{\substack{(g,q) \to (f,p) \\ q = P_{\mathcal{P}_n}(g)}} \frac{\langle \gamma, g-f \rangle - \langle \mu, q-p \rangle}{\|g-f\| + \|q-p\|}$$

$$\geq \limsup_{\substack{(g_\lambda,q) \to (f,p) \\ q = P_{\mathcal{P}_n}(g_\lambda)}} \frac{\langle \gamma, g_\lambda-f \rangle - \langle \mu, q-p \rangle}{\|g_\lambda-f\| + \|q-p\|}$$

$$= \limsup_{(f+\lambda, p+\lambda) \to (f,p)} \frac{\langle \gamma, f+\lambda-f \rangle - \langle \mu, p+\lambda-p \rangle}{\|f+\lambda-f\| + \|p+\lambda-p\|}$$

$$= \limsup_{\lambda \downarrow 0} \frac{\langle \gamma, \lambda \rangle - \langle \mu, \lambda \rangle}{2\lambda}$$

$$= \limsup_{\lambda \downarrow 0} \frac{\lambda(\gamma([0,1]) - \mu([0,1])}{2\lambda}$$

$$= \frac{\gamma([0,1]) - \mu([0,1]}{2} > 0.$$

By (5.13), this implies that

$$\mu([0,1]) < \gamma([0,1]) \implies \gamma \notin \widehat{D}^* P_{\mathcal{P}_n}(f,p)(\mu). \tag{5.14}$$

Similarly, to the proof of (5.14), by taking $g_\lambda = f - \lambda$, for $\lambda \downarrow 0$, in the limit (5.13), we can show

$$\mu([0,1]) > \gamma([0,1]) \implies \gamma \notin \widehat{D}^* P_{\mathcal{P}_n}(f,p)(\mu).$$

This proves part (i) of this theorem. Part (ii) follows from part (i) immediately.

Proof of (iii). Suppose $\langle \gamma, f - p \rangle < 0$. In the limit (5.13), we take a line segment direction as $g_\alpha = (1 - \alpha)f + \alpha p$, for $\alpha \downarrow 0$ with $\alpha < 1$. By part (iii) in Proposition 5.5, one has

$$p = P_{\mathcal{P}_n}((1 - \alpha)f + \alpha p), \text{ for any } \alpha \in (0, 1).$$

Since $(1 - \alpha)f + \alpha p \to f$, as $\alpha \downarrow 0$ with $\alpha < 1$, we have

$$\limsup_{\substack{(g,q)\to(f,p) \\ q=P_{\mathcal{P}_n}(g)}} \frac{\langle \gamma, g-f \rangle - \langle \mu, q-p \rangle}{\|g-f\| + \|q-p\|}$$

$$\geq \limsup_{\substack{(g_\alpha,q)\to(f,p) \\ q=P_{\mathcal{P}_n}(g_\alpha)}} \frac{\langle \gamma, g_\alpha-f \rangle - \langle \mu, q-p \rangle}{\|g_\alpha-f\| + \|q-p\|}$$

$$= \limsup_{\alpha\downarrow 0, \alpha<1} \frac{\langle \gamma, (1-\alpha)f + \alpha p - f \rangle - \langle \mu, p-p \rangle}{\|(1-\alpha)f + \alpha p - f\| + \|p-p\|}$$

$$= \limsup_{\alpha\downarrow 0, \alpha<1} \frac{-\langle \gamma, \alpha(f-p) \rangle}{\|\alpha(f-p)\|}$$

$$= \frac{-\langle \gamma, f-p \rangle}{\|f-p\|} > 0.$$

By (5.13), from this, it follows that

$$\langle \gamma, f - p \rangle < 0 \implies \gamma \notin \widehat{D}^* P_{\mathcal{P}_n}(f,p)(\mu), \text{ for any } \mu \in C^*. \qquad \square$$

Let $n$ be a nonnegative integer. By Theorem 5.4, the metric projection $P_{\mathcal{P}_n}: C[0, 1] \to \mathcal{P}_n$ is a single-valued mapping. So, it allows us to consider the Gâteaux directional differentiability of the metric projection operator $P_{\mathcal{P}_n}$.

**Proposition 5.11.** *Let n be a nonnegative integer. For any $f \in C[0, 1]$, we have*

$$P'_{\mathcal{P}_n}(f)(p) = p, \text{for any } p \in \mathcal{P}_n \setminus \{\theta\}.$$

*Proof.* For any $f \in C[0, 1]$ and for any $p \in \mathcal{P}_n \setminus \{\theta\}$, by the properties of the metric projection operator $P_{\mathcal{P}_n}$ in Proposition 5.5 and by $tp \in \mathcal{P}_n \setminus \{\theta\}$, for any $t > 0$, we have

$$P'_{\mathcal{P}_n}(f)(p) = \lim_{t\downarrow 0} \frac{P_{\mathcal{P}_n}(f+tp) - P_{\mathcal{P}_n}(f)}{t}$$

$$= \lim_{t\downarrow 0} \frac{P_{\mathcal{P}_n}(f) + tp - P_{\mathcal{P}_n}(f)}{t} = p. \qquad \square$$

**Acknowledgments.** The author is very grateful to Professor Boris S. Mordukhovich for his kind communications, valuable suggestions and enthusiasm encouragements in the development stage of this paper.

**Appendix**

**1**. **Some properties of the normalized duality mapping and generalized identity in general Banach spaces**

The inverse of the normalized duality mapping $J$ is denoted by $J^{-1}$, which is a set-valued mapping from $X^*$ to $X$ such that, for any $\varphi \in X^*$, we have

$$J^{-1}(\varphi) = \{x \in X : \varphi \in J(x)\}.$$

Then, for any $x \in J^{-1}(\varphi)$, we have $\|x\| = \|\varphi\|_*$. For $x, y \in X$, if $J(x) \cap J(x) \neq \emptyset$, then $x$ and $y$ are said to be generalized identical. For $x \in X$, the set of all its generalized identical points is denoted by $\Im(x)$. That is,

$$\Im(x) = \{y \in X : J(x) \cap J(x) \neq \emptyset\}.$$

The generalized identity has the following properties.

(a) $x \in \Im(x)$, for any $x \in X$;
(b) $\|y\| = \|x\|$, for any $y \in \Im(x)$.

We list some properties of the normalized duality mapping below for easy reference. For more details, one may see Sections 4.2–4.3 and Problem set 4.2 in [28] and [12, 17, 23, 25].

($J_1$). For any $x \in X$, $J(x)$ is nonempty, bounded, closed and convex subset of $X^*$;
($J_2$). $J$ is the identity operator in any Hilbert space $H$;
($J_3$). $J(\theta) = \theta^*$ and $J^{-1}(\theta^*) = \theta$;
($J_4$). For any $x \in X$ and any real number $\alpha$, $J(\alpha x) = \alpha J(x)$;
($J_5$). For any $x, y \in X$, $j(x) \in J(x)$ and $j(y) \in J(y)$, $\langle j(x) - j(y), x - y \rangle \geq 0$;
($J_6$). For any $x, y \in X$, $j(x) \in J(x)$, and $j(y) \in J(y)$, we have

$$2\langle j(y), x - y \rangle \leq \|x\|^2 - \|y\|^2 \leq 2\langle j(x), x - y \rangle;$$

($J_7$). $X$ is strictly convex if and only if $J$ is one-to-one, that is, for any $x, y \in X$,

$$x \neq y \implies J(x) \cap J(y) = \emptyset;$$

($J_8$). $X$ is strictly convex if and only if, for any $x, y \in X$ with $x \neq y$,

$$j(x) \in J(x) \text{ and } j(y) \in J(y) \implies \langle j(x) - j(y), x - y \rangle > 0;$$

($J_9$). $X$ is strictly convex if and only if, for any $x, y \in X$ with $\|x\| = \|y\| = 1$ and $x \neq y$,

$$\varphi \in J(x) \implies 1 - \langle \varphi, y \rangle > 0;$$

($J_{10}$). If $X^*$ is strictly convex, then $J$ is a single-valued mapping;
($J_{11}$). $X$ is smooth if and only if $J$ is a single-valued mapping;
($J_{12}$). If $J$ is a single-valued mapping, then $J$ is a norm to weak* continuous;
($J_{13}$). $X$ is reflexive if and only if $J$ is a mapping of $X$ onto $X^*$;
($J_{14}$). If $X$ is smooth, then $J$ is a continuous operator;
($J_{15}$). $J$ is uniformly continuous on each bounded set in uniformly smooth Banach spaces;
($J_{16}$) $X$ is strictly convex, if and only if $\Im(x) = x$, for every $x \in X$.

## 2. Proof of Proposition 2.1

**Proposition 2.1**. *Let $X$ be a Banach space and $C$ a nonempty closed and convex subset of $X$. For any $x \in X$ and $u \in C$, then, $u \in P_C(x)$, if and only if, there is $j(x - u) \in J(x - u)$ such that*

$$\langle j(x - u), u - y \rangle \geq 0, \quad \text{for all } y \in C. \tag{A.1}$$

*Proof.* Suppose that there is $j(x - u) \in J(x - u)$ such that (A.1) holds, for all $y \in C$. By the property ($J_6$) of the normalized duality mapping on Banach spaces, for any $j(x - u) \in J(x - u)$ and for any $y \in C$, by (A.1), we have

$$\begin{aligned}
\|x - y\|^2 - \|x - u\|^2 &\geq 2\langle j(x - u), (x - y) - (x - u) \rangle \\
&= 2\langle j(x - u), u - y \rangle \\
&\geq 0, \quad \text{for all } y \in C.
\end{aligned}$$

It implies that $u \in P_C(x)$. On the other hand, for $x \in X$, we consider the following two cases.

Case 1. $x \in X \backslash C$ and $u \in C$, if $u \in P_C(x)$, by Theorem 3.8.4 (it is for general normed spaces) in Zalinescu [30], there is $f \in X^*$ with $\|f\|_* = 1$ such that

$$\langle f, u - x \rangle = \|x - u\| \quad \text{and} \quad \langle f, y - u \rangle \geq 0, \quad \text{for all } y \in C. \tag{A.2}$$

Define $g \in X^*$ by $g = -\|x - u\| f$. From $\|f\|_* = 1$, we have

$$\|g\|_* = \|-\|x - z\| f\|_* = \|x - u\| \|g\|_* = \|x - u\|. \tag{A.3}$$

By the first equation in (A.2), it follows that

$$\langle g, x - u \rangle = \langle -\|x - u\| f, x - u \rangle = \|x - u\| \langle f, u - x \rangle = \|x - u\|^2. \tag{A.4}$$

Combining (A.3) and (A.4), we have $g \in J(x - u)$. By the second part of (A.2), we have

$$\langle g, u - y \rangle = \langle -\|x - z\| f, u - y \rangle = \|x - u\| \langle f, y - u \rangle \geq 0, \quad \text{for all } y \in C.$$

Then, taking $j(x - u) = g \in J(x - u)$, which satisfies (A.1).

Case 2. $x \in C$. In this case $P_C(x) = \{x\}$, which is a singleton. Then, for $u \in C$, we have that $u \in P_C(x)$ if and only if $u = x$. Then, $J(x - u) = J(\theta) = \{\theta^*\}$. It is clear (A.1) is satisfied. □

## 3. Proof of Proposition 3.1

**Proposition 3.1.** For any $a > 0$, we have

(a) $J^{-1}(\beta_a) = \Delta_a$;
(b) $J(y) = \beta_a$, for any $y \in \Delta_a^+$;
(c) $J(y) \supsetneq \{\beta_a\}$, for any $y \in \Delta_a \setminus \Delta_a^+$;
(d) $\beta_a = \bigcap_{x \in \Delta_a} J(x)$;
(e) $\Im(y) = \Delta_a$, for every $y \in \Delta_a^+$;
(f) $\Im(y) \supsetneq \Delta_a$, for every $y \in \Delta_a \setminus \Delta_a^+$.

*Proof.* Proof of (a). Let $y = (t_1, t_2, \ldots) \in l_1$. Suppose that $y \in J^{-1}(\beta_a)$. It follows that $\beta_a \in J(y)$. This implies

$$\|y\|_1^2 = \langle \beta_a, y \rangle = \|\beta_a\|_\infty^2 = a^2.$$

Since $\langle \beta_a, y \rangle = a \sum_{n=1}^\infty t_n = a^2$, it yields that $\sum_{n=1}^\infty t_n = a$. Then, it follows that $\|y\|_1 = \sum_{n=1}^\infty |t_n| = a = \sum_{n=1}^\infty t_n$, it implies that $t_n \geq 0$, for $n = 1, 2, \ldots$. That is, $y \in K_1$ and $\|y\|_1 = a$, which implies $y \in \Delta_a$.

On the other hand, if $y = (t_1, t_2, \ldots) \in \Delta_a$ with $y \in K_1$ and $\|y\|_1 = a$. Then, $\|y\|_1^2 = a^2 = \langle \beta_a, y \rangle = \|\beta_a\|_\infty^2$. This implies $\beta_a \in J(y)$, which reduces $y \in J^{-1}(\beta_a)$. So, (a) is proved.

Proof of (b). Let $y = (t_1, t_2, \ldots) \in \Delta_a^+$. Then, $t_n > 0$, for all $n$ and $\|y\|_1 = a$. For any $j(y) = (v_1, v_2, \ldots) \in J(y)$, we have $\|j(y)\|_\infty^2 = \langle j(y), y \rangle = \|y\|_1^2 = a^2$. It follows that

$$a = \|j(y)\|_\infty = \sup_{1 \leq n < \infty} |v_n|.$$

Then, from $\langle j(y), y \rangle = \sum_{n=1}^\infty t_n v_n = a^2$, $t_n > 0$, for all $n$, and $\sum_{n=1}^\infty t_n = a$, it implies $v_n = a$, for all $n = 1, 2, \ldots$. This is $j(y) = \beta_a$.

Proof of (c). Let $y = (t_1, t_2, \ldots) \in \Delta_a \setminus \Delta_a^+$. There is a positive integer $m$ such that $t_m = 0$. We define $\gamma = (u_1, u_2, \ldots) \in l_\infty$ satisfying

$$u_n = \begin{cases} 0, & \text{if } n = m, \\ a, & \text{if } n \neq m, \end{cases} \quad \text{for all } n.$$

One sees that $\|\gamma\|_\infty^2 = \langle \gamma, y \rangle = a^2 = \|y\|_1^2$. This implies $\gamma \in J(y)$ with $\gamma \neq \beta_a$, which proves (c).

Proof of (d). By part (a), we have $\beta_a \in \bigcap_{x \in \Delta_a} J(x)$. On the other hand, for any $\varphi = (u_1, u_2, \ldots) \in l_\infty$ with $\varphi \neq \beta_a$, then, there is a positive integer $m$ such that $u_m \neq a$. We define $y = (t_1, t_2, \ldots) \in \Delta_a$ satisfying

$$t_n = \begin{cases} a, & \text{if } n = m, \\ 0, & \text{if } n \neq m, \end{cases} \quad \text{for all } n.$$

One sees that $\langle \varphi, y \rangle = au_m \neq a^2 = \|y\|_1^2$. This implies $\varphi \notin J(y)$, which proves (d).

Proof of (e). By part (b), for $y \in \Delta_a^+$, we have $J(y) = \beta_a$. Then, by part (a), we have $J^{-1}(\beta_a) = \Delta_a$, this implies $\mathfrak{J}(y) = \Delta_a$, for every $y \in \Delta_a^+$.

Proof of (f). Let $y = (t_1, t_2, \ldots) \in \Delta_a \setminus \Delta_a^+$. By part (d), we have $\mathfrak{J}(y) \supseteq \Delta_a$, for every $y \in \Delta_a \setminus \Delta_a^+$. There is a positive integer $m$ such that $t_m = 0$. We define $\gamma = (u_1, u_2, \ldots) \in l_\infty$ satisfying

$$u_n = \begin{cases} -a, & \text{if } n = m, \\ a, & \text{if } n \neq m, \end{cases} \text{ for all } n.$$

One sees that $\|\gamma\|_\infty^2 = \langle \gamma, y \rangle = a^2 = \|y\|_1^2$. This implies $\gamma \in J(y)$ with $\gamma \neq \beta_a$. Let $x = (s_1, s_2, \ldots)$ defined by

$$s_n = \begin{cases} -\frac{1}{2}, & \text{if } n = m, \\ \frac{t_n}{2}, & \text{if } n \neq m, \end{cases} \text{ for all } n.$$

One has that $\|\gamma\|_\infty^2 = \langle \gamma, x \rangle = a^2 = \|x\|_1^2$. This implies $\gamma \in J(x) \cap J(y)$ with $\gamma \neq \beta_a$. This implies $x \in \mathfrak{J}(y)$. However, $x \notin \Delta_a$. This proves part (f). □

### 4. A Proof of Theorem 5.4

**Theorem 5.4.** *Let $n$ be a nonnegative integer. The metric projection $P_{\mathcal{P}_n}: C[0, 1] \to \mathcal{P}_n$ is a single-valued mapping.*

*Proof.* For any $f \in C[0, 1]$, by Proposition 5.2, we have $P_{\mathcal{P}_n}(f) \neq \emptyset$. So, in this theorem, we only need to show that $P_{\mathcal{P}_n}(f)$ is a singleton by using the Chebyshev's Equioscillation Theorem.

It is clear that, for any $p \in \mathcal{P}_n$, we have $P_{\mathcal{P}_n}(p) = p$, which shows that $P_{\mathcal{P}_n}$ is a single-valued mapping on $\mathcal{P}_n$. Let $f \in C[0, 1] \setminus \mathcal{P}_n$. Assume there are $p, q \in P_{\mathcal{P}_n}(f)$. Then, $\|f - p\| = \|f - q\|$. Take an $n$-Chebyshev set of $f$ with respect to $p$ as $S(f, p) = \{t_0, t_1, \ldots, t_{n+1}\}$. Without loose of the generality, we suppose $\epsilon = 1$. By the Chebyshev's Equioscillation Theorem, this implies

$$f(t_i) - p(t_i) = (-1)^i \|f - p\|, \ i = 0, 1, 2, \ldots, n+1.$$

By $\|f - p\| = \|f - q\|$, this implies that, for any $i$ in $\{0, 1, 2, \ldots, n+1\}$, we have

$$f(t_i) - q(t_i) \leq \|f - p\|, \text{ if } i \text{ is even},$$

and
$$f(t_i) - q(t_i) \geq -\|f - p\|, \text{ if } i \text{ is odd}.$$

This implies

$$q(t_i) - p(t_i) \geq 0, \text{ if } i \text{ is even},$$

and
$$q(t_i) - p(t_i) \leq 0, \text{ if } i \text{ is odd}.$$

Then, for any $i = 0, 1, 2, \ldots, n-1$, there is $s_i \in [t_i, t_{i+2}]$ such that

$$(q - p)'(s_i) = 0, \text{ for any } i = 0, 1, 2, \ldots, n.$$

At the points $s_0, s_1, \ldots, s_n$, the polynormal $q - p$ takes local minimum and local maximum alternatively. Then, we prove that $q - p$ is a constant function by the following two cases.

Case 1. There is $i \in \{0, 1, 2, \ldots, n-2\}$ such that $s_i = s_{i+1}$. Since $q - p$ takes local minimum at one of these two points $s_i$ and $s_{i+1}$ (on $[t_i, t_{i+2}]$ or $[t_{i+1}, t_{i+3}]$, respectively), and $q - p$ takes local maximum at other one, it implies that $q - p$ must be constant in $[t_{i+1}, t_{i+2}]$. Since $t_{i+1} < t_{i+2}$, this implies that $q - p$ is a constant function.

Case 2. $s_0 < s_1 <, \ldots < s_{n-1}$. In this case, the polynormal $(q - p)'$ has $n$ distinct zeros in $[0, 1]$. Since $(q - p)'$ is a polynormal with degree less than or equal to $n-1$, this implies that $(q - p)' = 0$. Then, it follows that $q - p$ is a constant function.

Suppose $q - p = \lambda$, for some real number $\lambda$. That is, $q = p + \lambda$. By $\|f - p\| = \|f - q\|$, we have

$$\|f - p\| = \|f - q\| = \|f - (p + \lambda)\| = \|f - p - \lambda\|.$$

This implies $\lambda = 0$. Hence, $q = p$. □